\documentclass{paper} 

\usepackage{amssymb}
\usepackage{amsmath}
\usepackage{amsthm}
\usepackage{color}
\usepackage{graphicx}
\usepackage{enumerate}
\usepackage{subfigure}

\newtheorem{thm}                   {Theorem}    
\newtheorem{prop}         [thm]{Proposition}
\newtheorem{lemma}      [thm]{Lemma}
\newtheorem{cor}           [thm]{Corollary}

\theoremstyle{definition} 
\theoremstyle{definition} \newtheorem{remark}[thm]{Remark}
\theoremstyle{definition} \newtheorem{problem}[thm]{Problem}

\begin{document}

\vbox to 50pt {}

\centerline{
{\Large  
{\bf Kronecker covers,} {\it V}{\bf -construction, 
}}}
\medskip
\centerline{
{\Large 
{\bf unit-distance graphs and isometric}}}
\medskip
\centerline{
{\Large 
{\bf 
point-circle configurations\footnote
{The authors acknowledge partial funding of this research via ARSS of Slovenia, grants: P1-0294 and N1-0011: GReGAS, supported 
in part by the European Science Foundation.}
}}}
\bigskip \bigskip

\centerline{
G\'abor G\'evay 
}
\centerline{
Bolyai Institute, University of Szeged
}
\centerline{
Aradi v\'ertan\'uk tere 1, H-6720 Szeged
}
\centerline{
Hungary
}
\centerline{
{\tt gevay@math.u-szeged.hu}
}

\smallskip

\centerline{and}

\smallskip

\centerline{ 
Toma\v z Pisanski
}
\centerline{ 
Faculty of Mathematics and Physics
}
\centerline{ 
University of Ljubljana, Jadranska 19, 1111 Ljubljana
}
\centerline{ 
Slovenia
}
\centerline{
{\tt Tomaz.Pisanski@fmf.uni-lj.si}
}
\medskip

\date{\today} 

\bigskip
\medskip
\centerline{
{\bf Abstract}
}

\medskip
\noindent
{\small We call a polytope $P$ of dimension 3 \emph{admissible} if it has the following two properties: (1) for each vertex of $P$ the set of 
its first-neighbours is coplanar; (2) all planes determined by the first-neighbours are distinct. It is shown that the Levi graph of a point-plane 
configuration obtained by $V\!$-construction from an admissible polytope $P$ is the Kronecker cover of its 1-skeleton. We investigate the 
combinatorial nature of the $V$-construction and use it on unit-distance graphs to construct novel isometric point-circle configurations. In 
particular, we present an infinite series whose all members are subconfigurations of the renowned
Clifford configurations.}
\medskip

\noindent
{\small 
\medskip

\noindent
\emph{Math.\ Subj.\ Class.: 05B30, 51A20, 52B10, 52C30}}

\bigskip

\section{Introduction.}
In this paper we investigate and carry over from polytopes to graphs the so-called $V\!$-construction, which was originally introduced in 
\cite{GG09}.  In this process we explain the construction in terms of the canonical double cover, also called the Kronecker cover of graphs. 
The reader is referred to \cite{IP08} for graph coverings and to the monographs \cite{Gru09,PS} for the background on configurations and 
their Levi graphs.

The first author used convex 3-polytopes in order to define a construction of geometric point-plane configurations in the following way~\cite{GG09}.

Let $P$ be a polytope of dimension 3 with the property that for each vertex $v$ the set of its first neighbours $P(v)$ is coplanar. 
In particular, this will always be true if the graph (or 1-skeleton) of the polytope $P$ is trivalent. There are several other classes of 
polytopes that have this property. Furthermore, we assume that all planes obtained in this way are distinct. In particular, this condition 
rules out bipyramids such as the octahedron. Let us call such a polytope \emph{admissible}. 

\begin{prop}
Each 3-polytope with trivalent 1-skeleton is admissible.
\end{prop}

\noindent
\emph{Proof.} 
Since each vertex of a 3-polytope with trivalent 1-skeleton has exactly three first-neighbours, they are clearly coplanar. An easy argument shows 
that if two trivalent vertices of a 3-polytope $P$ share the same set of first-neighbours, then the 1-skeleton of $P$ itself cannot be trivalent. 
\hfill $\square$
\bigskip

Let $S(P)$ denote the set of such planes as above, if they exist. If $V(P)$ denotes the set of vertices of $P$, then the pair $((V (P), S(P))$ defines a 
geometric incidence structure of points and planes with the usual incidence. We call this procedure the \emph{geometric $V\!$-construction}. If the 
1-skeleton of $P$ is a regular, say $k$-valent graph, then each point of the configuration will sit on $k$ planes.  It immediately follows from the 
definition that each plane contains exactly $k$ points. Let $n$ be the number of vertices of $P$. Therefore, combinatorially, the incidence structure 
is an $(n_k)$ configuration. 

A natural question is that what is the Levi graph of such a configuration. Recall that the \emph{Levi graph} $L(C)$ of a configuration $C$ is 
a bipartite graph whose bipartition classes consist of the points and ``blocks" of $C$, respectively, and two points in $L(C)$ are adjacent if 
and only if the corresponding point and ``block" in $C$ are incident. Levi graphs are useful tools in studying configurations, because of the
following property~\cite{Cox50}.

\begin{lemma}
A configuration $C$ is uniquely determined by its Levi graph $L(C)$. \label{Levi}
\end{lemma}

Another, much more difficult question is, if we can find any conditions under which such a combinatorial configuration may be realized 
geometrically as configuration of points and lines. On the other hand, it may happen that a configuration can be realized in both a point-line 
and a point-plane version (cf.\ our Example in  Section~\ref{comb_V_con}).  In Section~\ref{PointCircle} and~\ref{KneserConfig} we also 
present examples for configurations of which both point-line and point-circle realization exist. 

Point-circle configurations themselves are also interesting, since, in contrast to the point-line configurations, relatively little is known about them.
The most notable achievement in this respect is undoubtedly Clifford's infinite series of configurations, going back to 1871~\cite{Cox61,Gru09}.
In the last two sections we present a new construction of Clifford's configurations, as well as three new infinite series of point-circle configurations.

We note that the construction introduced in \cite{GG09} is more general than needed here: instead of 3-dimensional
polytopes one can take $d$-dimensional polytopes and accordingly, instead of planes one should consider hyperplanes.
Also, instead of first-neighbours it is possible to consider second-neighbours. However, we do not consider these 
aspects of the $V\!$-construction here.

\section{Combinatorial $V\!$-construction.} \label{comb_V_con}

Let us generalize and carry out the $V\!$-construction on the abstract level. 

To any regular graph $G$ we may associate a combinatorial configuration. For a vertex $v$ of $G$, denote by $N(v)$ the set of vertices adjacent
to $v$. Then take the family $S(G)$ of these vertex-neighbourhoods:
$$
S(G) = \{N(v)\, | \, v \in V (G) \}.
$$
The triple $(V (G), S(G), \in)$ defines a combinatorial incidence structure, underlying the geometric configuration of points and planes 
for any 3-polytope $P$ whose 1-skeleton is $G$. We shall denote this structure by 
$N(G)$.

We note that a closely related construction occurs in the context of combinatorial geometries~\cite{LPL,Van}. 

The following general result establishes a connection between Levi graphs and Kronecker covers. It will play a central role in our constructions 
presented in the rest of the paper. First, we recall that a graph $\widetilde G$ is said to be the \emph {Kronecker cover} (or  \emph {canonical 
double cover}) of the graph $G$ if there exists a $2:1$ surjective homomorphism $f : \widetilde G \rightarrow G$ such that for every vertex 
$v$ of $\widetilde G$ the set of edges incident with $v$ is mapped bijectively onto the set of edges incident with $f (v)$~\cite{IP08}. 

\begin{thm} \label{MainThm}
Let $G$ be a graph on $n$ vertices and let $L$ be the Levi graph of the incidence structure $N(G)$. If no two vertices of $G$ have the same 
neighbourhood, then $L$ is the Kronecker cover of $G$. 
\end{thm}

\noindent
\emph{Proof.}   Under the assumption that no two vertices have the same set of neighbours, all sets $N(v)$, for $v \in V(G)$, are distinct.  
Therefore the set of vertices of $L$ consists of $V$ and $\{N(v)| v \in V(G)\}$. Each edge $e = uv$ from $G$ gives rise to two edges: 
$uN(v)$ and $vN(u)$. Hence $L$ is a Kronecker cover of $G$. If $|V(G)| \neq |\{N(v)| v \in V(G)\}|$, the argument fails. 
\hfill $\square$
\bigskip

\noindent
Some direct consequences of Theorem \ref{MainThm} for Levi graphs are as follows.

\begin{prop}
Let $G$ be a graph on $n$ vertices and let $L$ be the Levi graph of the incidence structure $N(G)$. The graph $L$ is connected if and only if 
$G$ is connected and non-bipartite. 
\end{prop}

\noindent
\emph{Proof.}  In case the graph $G$ has no two vertices with a common neighborhood, the result follows from a well-known property 
of the Kronecker cover, see Proposition 1 of \cite{IP08}. If this is not the case, the construction of $L$ may be performed in two steps. 
First we construct the Kronecker cover over $G$ and then identify some pairs of vertices, such as $N(u)$ and $N(v)$, in case $N(u) = N(v)$.  
Such an identification may occur if and only if the vertices $u$ in $v$ are in the same bipartition set. This means that in the Kronecker 
cover only vertices in the same connected component may be identified. 
\hfill $\square$
\smallskip

\begin{prop}
Let $G$ be a regular graph of valency $k$ on $n$ vertices and let $L$ be the Levi graph of the incidence structure  $N(G)$. 
Then $N(G)$ is an abstract $(n_k)$ point-line configuration if and only if $L$ contains no cycle of length 4. 
\end{prop}

\noindent
\emph{Proof.}  In the Kronecker cover odd cycles of length $r$ lift to cycles of length $2r$, while even cycles lift to two
cycles of the same length. Hence the girth of the Kronecker cover is 4 if and only if the original graph contains a 4-cycle. 
Since Kronecker cover is bipartite, the alternative means girth at least 6. \hfill $\square$
\bigskip

By analogy with geometric $V\!$-construction, we call a graph $G$ \emph{admissible} if no two of its vertices
have a common neighborhood. Recall that a configuration is \emph{combinatorially self-polar} if there exists an automorphism
of order two of its Levi graph interchanging the two parts of bipartition; see for instance \cite{PS}.

\begin{thm}
A configuration that is obtained by $V\!$-construction from an admissible graph $G$ is
combinatorially self-polar.
\end{thm}

\noindent
\emph{Proof.}  By our previous discussion the Levi graph of this configuration is a Kronecker cover over $G$. 
The involution that switches at the same time the vertices in each fiber is self-polarity. This follows from the fact
that any double cover is a regular cover.  \hfill $\square$
\bigskip

\noindent
We shall use the following result, which is an easy consequence of Proposition 1 in~\cite{IP08} and our Theorem~\ref{MainThm}. 

\begin{prop} \label{bipartite}
Let $C$ be a configuration obtained from a graph $G$ by $V\!$-con\-struc\-tion. Then the Levi graph $L(C)$ is bipartite. 
If $G$ is bipartite, then $L(C)$ consists of two disjoint copies of $G$.
\end{prop} 

\begin{cor} \label{TwoCopies}
Applying the $V\!$-construction to the Levi graph of configuration $C$ from Proposition \ref{bipartite} results in a configuration $C'$ 
which consists of two disjoint copies of $C$.
\end{cor}

\begin{figure}[h!] 
  \begin{center}
   \includegraphics[width=0.66\textwidth]{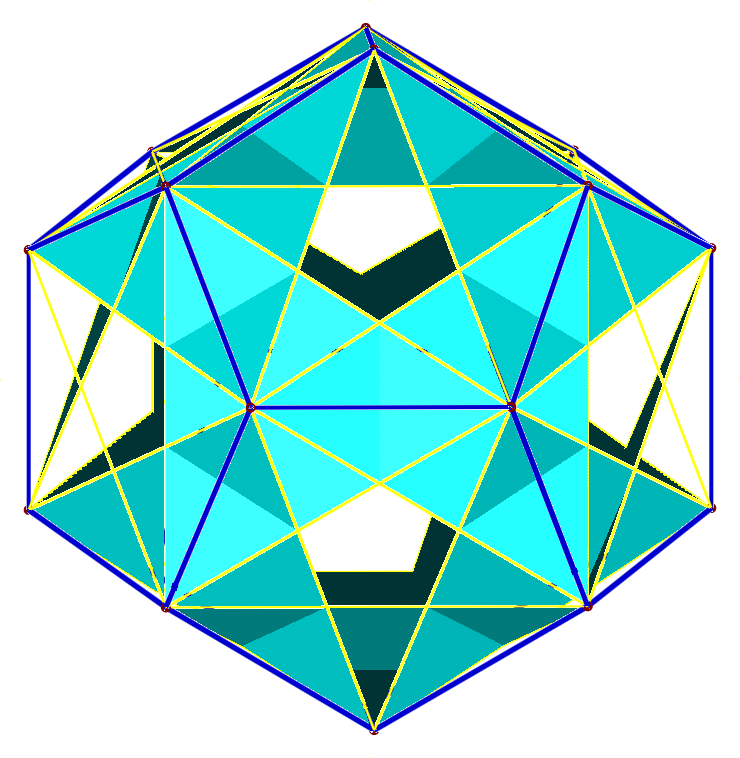}
  \end{center}
  \caption{The $(20_3)$ point-plane configuration obtained from the regular dodecahedron by the $V\!$-construction.}
  \label{Plane(20_3)}  
\end{figure}

We conclude this section with the following example. Let $G$ be the the dodecahedron graph. Then $N(G)$ is a configuration $(20_3)$. 
If $G$ is embedded in $\mathbb E^3$ as the 1-skeleton of the regular dodecahedron, then $N(G)$ is realized as a geometric point-plane 
configuration (see Figure~\ref{Plane(20_3)}). 

We note that the same configuration is obtained by taking the 20 vertices and the planes spanned by the 20 triangular faces of either the 
\emph {small ditrigonal icosidodecahedron} or the \emph {great ditrigonal icosidodecahedron} (these polyhedra belong to the class of the 
53 non-regular non-convex uniform polyhedra~\cite{CLM, Har}). 

\begin{figure}[h!] 
 \begin{center}
 \subfigure[]{
  \includegraphics[width=0.45\textwidth]{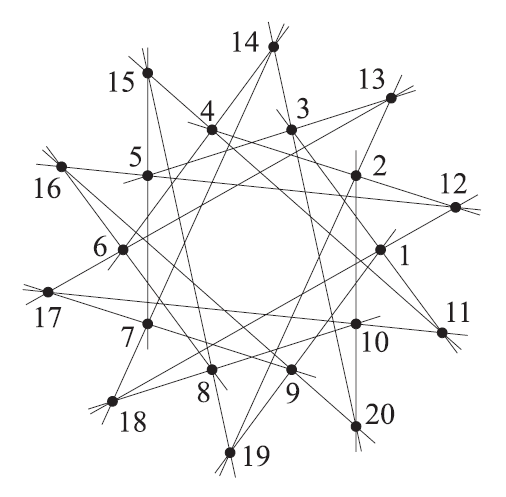}}
         \hskip 20pt
  \subfigure[]{\hskip -4pt
  \includegraphics[width=0.455\textwidth]{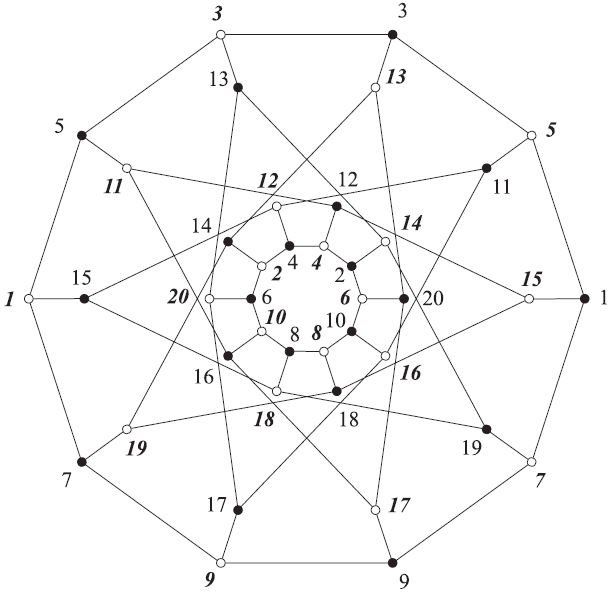}}
 \end{center}
\caption{The flag-transitive triangle-free point-line configuration $(20_3)$ (a),
and its Levi graph (b).}
\label{Line(20_3)}
\end{figure}

On the other hand, we know that the Kronecker cover of the dodecahedron graph is isomorphic with the Levi graph of the $(20_3)$ point-line 
configuration which is unique with the properties that it is triangle-free and flag-transitive~\cite{BGPZ} (Figure~\ref{Line(20_3)})  (see also 
Figure 1 in~\cite{BBP}). Thus we can see that $N(G)$ can be realized geometrically as both a point-plane and a point-line configuration.  
As we show in the next section, a realization as a point-circle configuration may be of interest. 
\section{$V\!$-construction and configurations of points and circles.} \label{PointCircle}

In~\cite{GG09} it was observed that certain point-plane configurations obtained from a 3-polytope $P$ by the $V\!$-construction could also 
be realized by points and circles. A simple necessary condition for this is that for each vertex $v$ of $P$, the set of the first-neighbours of $v$ 
forms a concyclic set, i.e.\ one can draw a circle through its points. Moreover, such point-circle configurations can be carried over to the plane, 
using stereographic projection. Here the well-known property is used that the stereographic projection is a circle-preserving map, see for 
instance~\cite{HCV} (also~\cite{Cox61, Hah}).

\begin{lemma}
Under stereographic projection from the sphere $S$ to the plane $\Sigma$ the image of any circle on $S$
is a circle on $\Sigma$. 
\end{lemma} 

\noindent
A straightforward application of this leads to the following result.

\begin{thm}
Any point-circle configuration on the sphere gives rise to a planar point-circle configuration.
\end{thm}

\noindent
Here we explicitly state the result that is presented already in \cite{GG09} (see Table 1 and 2 there), and follows readily from the above Theorem.

\begin{cor}\label{Platonic}
The $V\!$-construction of any Platonic or Archimedean polyhedron except for the octahedron
gives rise to a planar point-circle configuration.
\end{cor} 

\begin{figure}[h!] 
  \begin{center}
   \includegraphics[width=0.66\textwidth]{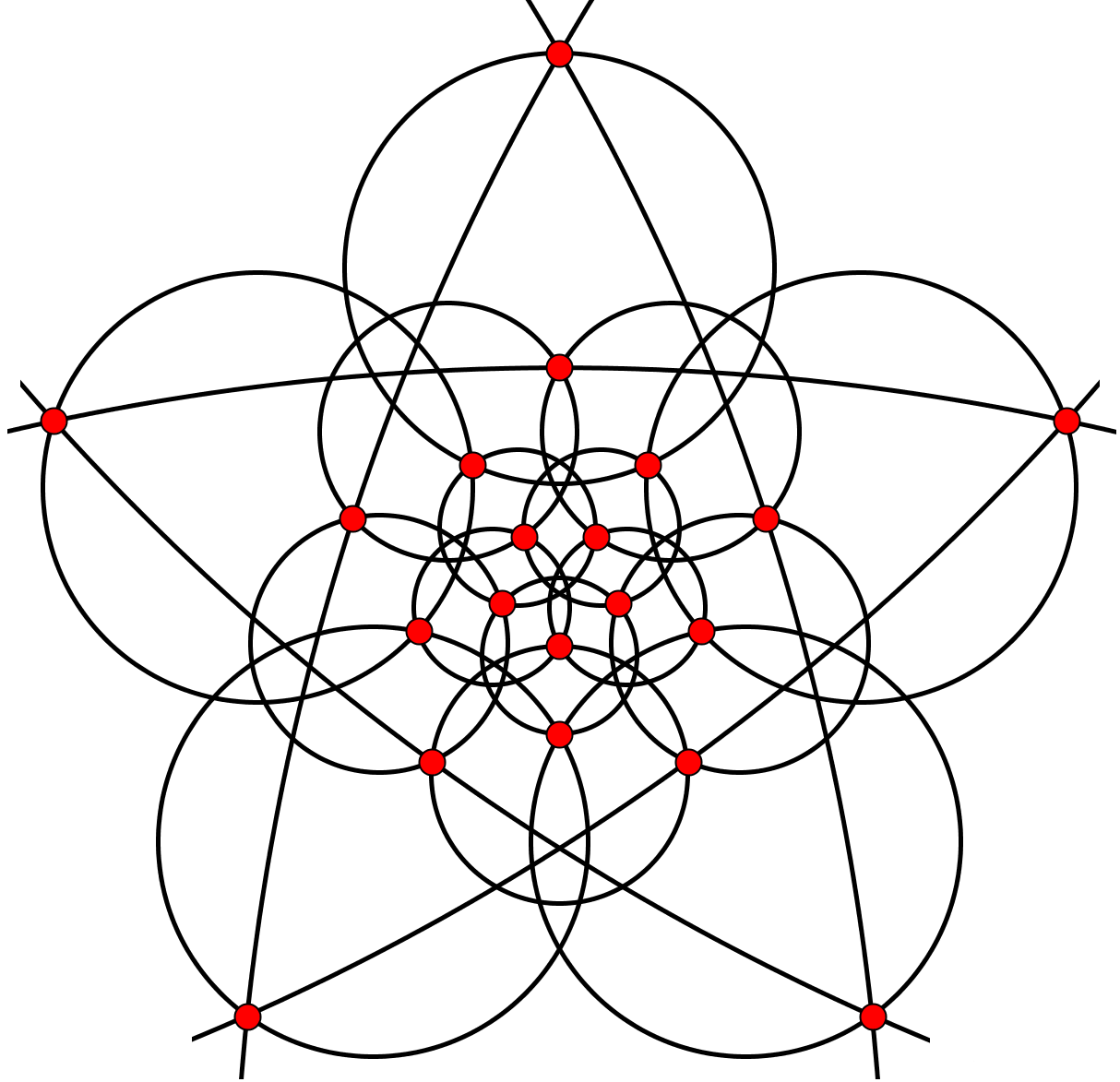}
  \end{center}
  \caption{$(20_3)$ point-circle configuration obtained from the regular dodecahedron by $V\!$-construction.}
  \label{V(Dod)}  
\end{figure}

An example obtained from the regular dodecahedron is depicted in Figure~\ref{V(Dod)}.  (Note that together with this, we have three distinct 
geometric realizations of one and the same abstract configuration of type $(20_3)$; cf.\ Figures~\ref{Plane(20_3)} and~\ref{Line(20_3)}.)

We remark that applying highly symmetric polytopes as a ``scaffolding" for the construction of spatial point-line configurations 
is extensively used in~\cite{GG13}.

In what follows we consider some other cases of $V\!$-construction that also give rise to planar point-circle configurations.

\begin{prop} \label{n_3}
Any $(n_3)$ configuration can be realized by points and circles in the plane.
\end{prop}

\noindent
\emph{Proof.} We may place the $n$ points in the plane in general position, in such a way that no three lie on a line and no four lie on a circle. 
Obviously combinatorial lines can be realized by circles, and the combinatorial incidence is carried over to a geometric point-circle incidence. 
\hfill $\square$
\bigskip

The $(n_3)$ point-circle configurations have an important property that is not shared by all point-circle configurations; namely, they are 
movable. To see this notion, we should consider that in the simplest case our point-circle configurations are constructed in the Euclidean 
plane $\mathbb E^2$. However, by adding to $\mathbb E^2$ a single point at infinity, we may consider them as lying in the \emph{inversive 
plane}~\cite{Cox61}. In this latter case, we say that a point-circle configuration is \emph{rigid} if its geometric realizations form a single 
class under circle-preserving transformations. 

We note that point-circle configurations can also be considered on the extended complex plane; in this case the circle-preserving transformations 
are just the M\"obius transformations, i.e.\ fractional linear transformations~\cite{Hah}. Incidentally, they play an important role in the so-called 
Lombardi drawings of graphs, an idea not totally unrelated to point-circle configurations and studied by D.\ Eppstein and his co-workers, for 
instance in~\cite{DEGKN}. 

A configuration that is not rigid is called \emph{movable} (cf.\ the notion of movability of point-line configurations, as defined in~\cite{Gru09}). 
Having defined this notion, the following statement is straightforward.

\begin{prop}
Any $(n_3)$ point-circle configuration is movable.
\end{prop}

\noindent
We note that movability is not a general property even for point-line configurations; 
for example, some classes of movable $(n_4)$ configurations were discovered just 
recently~\cite{Ber, BBGP}.

There is another property that distinguishes  $(n_3)$ point-circle configurations among all configurations. In general, the circles may be 
of different size. Let $r$ be the number of radii used in this construction. If $r = 1$, all circles are of the same size, and the configuration 
is called an \emph {isometric point-circle configuration}. It is not clear which $(n_3)$ configurations can be realized as planar isometric 
point-circle configurations. 

There is a large class of graphs that yields by V-construction isometric point-circle configurations in a natural way. 
These are the \emph {unit-distance graphs}, i.e.\ graphs whose all edges have the same length (cf.~\cite{HP10, HP12, ZHP}).

\begin{thm}
Let $G$ be a regular $d$-valent graph that is a unit-distance graph on $n$ vertices in the plane. 
Then $S(G)$  is an $(n_d)$ configuration, realizable as an isometric point-circle configuration.
\label{thm:ud}
\end{thm}

\noindent
\emph{Proof.} The points of the configuration are the vertices of the graph, as drawn in the plane. Unit-distance property implies 
that for each vertex, the set of its first-neighbours forms a concyclic set; furthermore, all these circles are of the same size. 
\hfill $\square$
\medskip

\begin{figure}[h!] 
  \begin{center}\hskip 4pt
  \subfigure[]{\hskip -4pt 
  \includegraphics[width=0.395\textwidth]{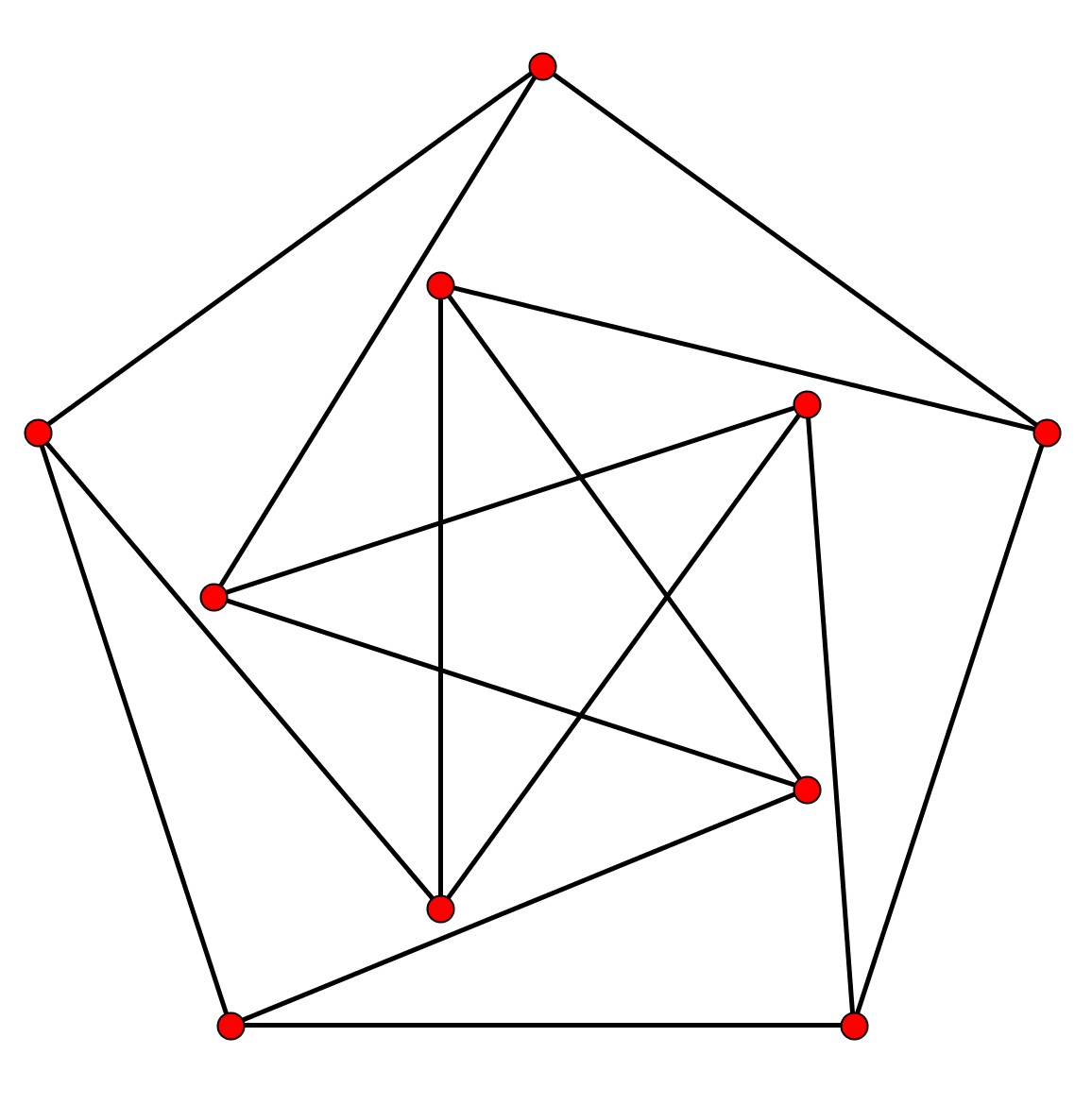}} \hskip 20pt
 \subfigure[]{\hskip -4pt 
   \includegraphics[width=0.415\textwidth]{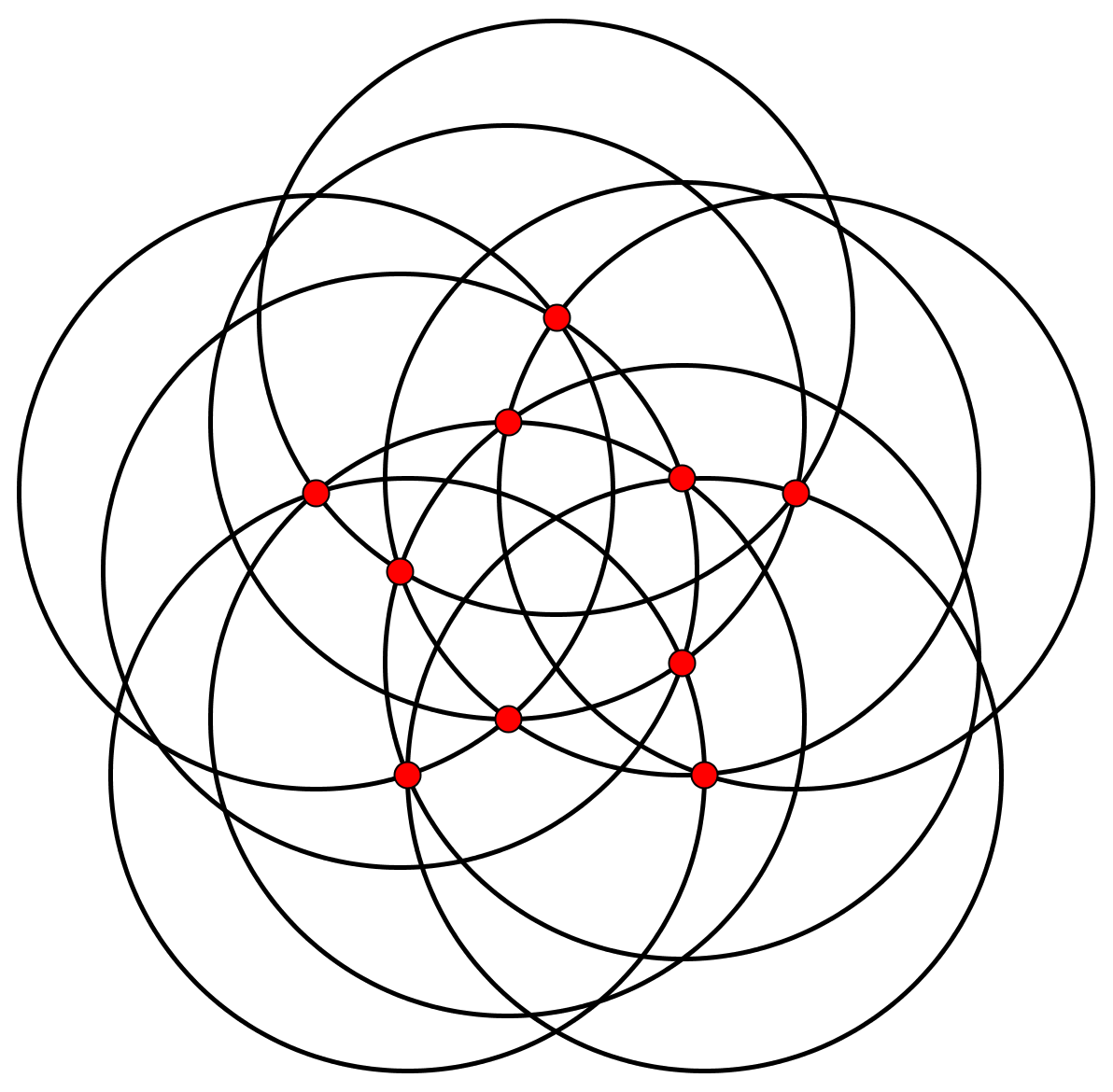}}
  \end{center}
  \caption{Applying the $V\!$-construction to a unit-distance representation of the Petersen graph (a) yields an isometric point-circle realization 
                 of Desargues' configuration (b).}
  \label{UD-Desargues}  
\end{figure}
An interesting example is as follows. We know unit-distance representations of the Petersen graph~\cite{HP12, ZHP}; on the other hand, it is 
well known that the Kronecker cover of the Petersen graph is the Desargues graph~\cite{IP08} (which, in turn, is the Levi graph of the Desargues 
configuration~\cite{Cox50}). Thus, on account of Theorem~\ref{MainThm}, the $V\!$-construction on a unit-distance representation of the Petersen 
graph yields an isometric point-circle realization of the Desargues configuration (see Figure~\ref{UD-Desargues}).

\begin{figure}[h] 
  \begin{center}\hskip 4pt
  \subfigure[]{\hskip -4pt 
  \includegraphics[width=0.425\textwidth]{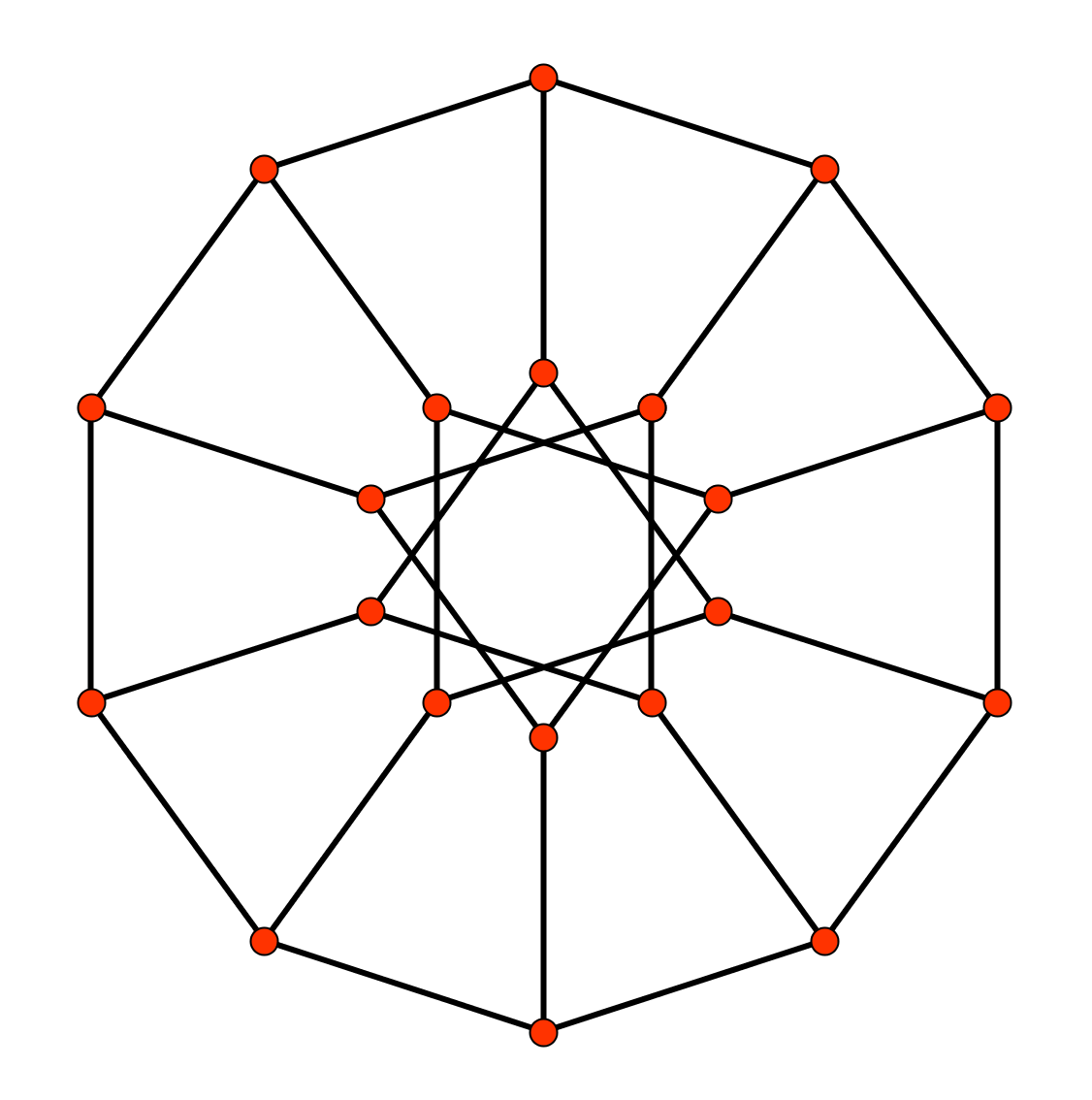}} \hskip 20pt
 \subfigure[]{\hskip -4pt 
   \includegraphics[width=0.433\textwidth]{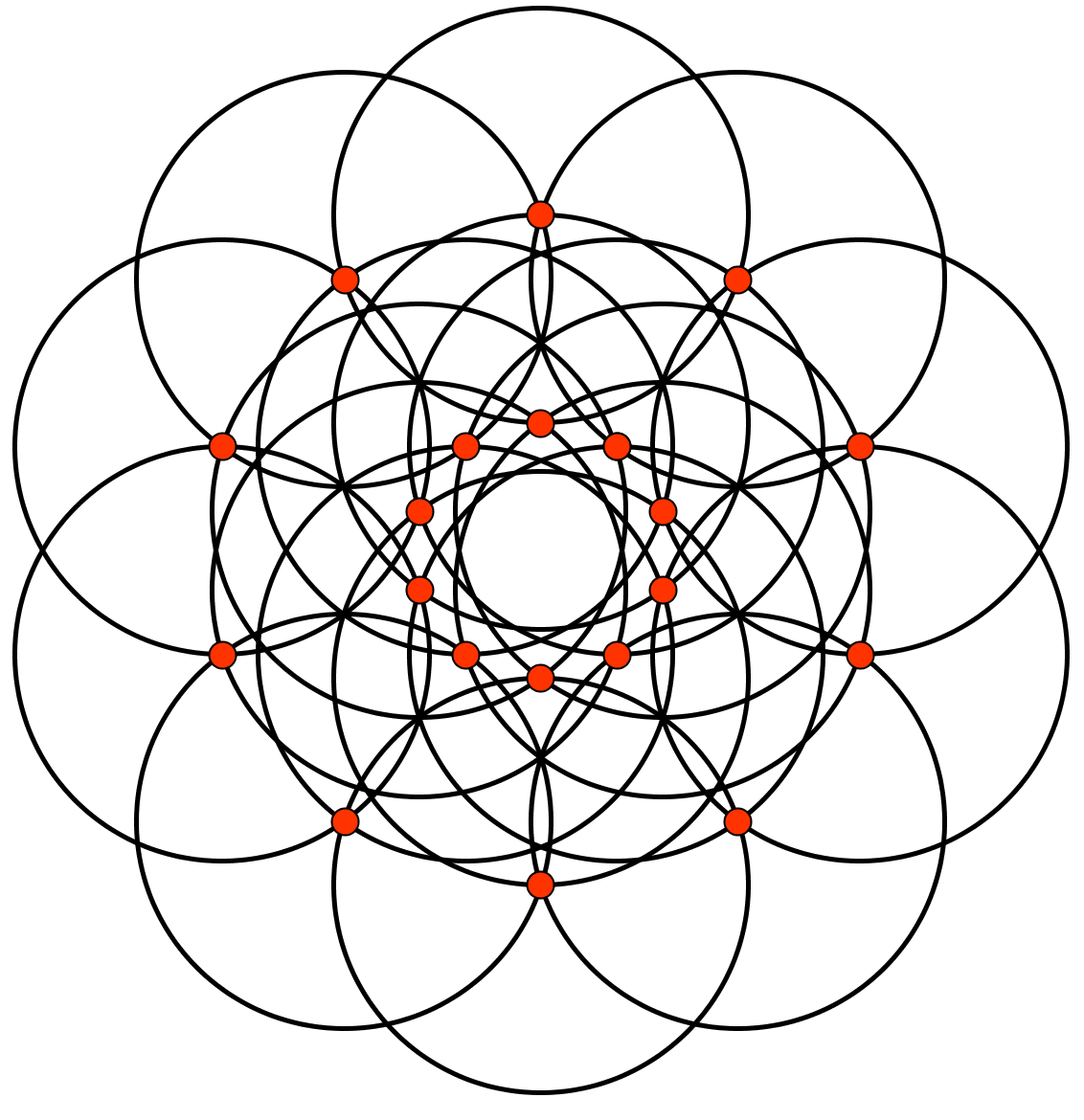}}
  \end{center}
  \caption{Applying the $V\!$-construction to a unit-distance representation of the Desargues graph (a) yields a configuration 
                 consisting of two disjoint copies of another isometric realization of the point-circle Desargues configuration (b). 
                 One of the two copies may be viewed in Figure \ref{NUD-Desargues}.}
  \label{DoubleDesargues}  
\end{figure}

We remark that the Desargues graph also has a unit-distance representation~\cite{ZHP}. Thus one may also apply to it the $V\!$-construction, 
so as to obtain an isometric point-circle configuration. By our Corollary~\ref{TwoCopies}, this $(20_3)$ configuration decomposes into two 
disjoint copies of the $(10_3)$ Desargues point-circle configuration (see Figure~\ref{DoubleDesargues}, where the construction yields the 
two $(10_3)$ copies in centrally symmetric position with respect to their common centre).

\section{$V\!$-construction on $d$-cubes and Clifford's point-circle con\-fig\-urations} \label{Clifford}

\noindent
The particular case of the cube in Corollary~\ref{Platonic} can be extended to the whole class of $d$-cubes. 
Because of its interesting connections, we discuss here the general case in some detail.

We recall that the infinite series of Clifford's point-circle configurations is associated to his renowned chain of theorems.
By Coxeter~\cite{Cox50}, these theorems can be formulated as follows.

\begin{thm}[{\bf Clifford's chain of theorems}]
{$\phantom{.}$}\newline
\vskip -11pt
{\rm (1)} Let $\sigma_1$, $\sigma_2$, $\sigma_3$, $\sigma_4$ be four circles of general position through a point $S$. Let $S_{ij}$ 
be the second intersection of the circles $\sigma_i$ and $\sigma_j$. Let  $\sigma_{ijk}$ denote the circle $S_{ij}S_{ik}S_{jk}$. 
Then the four  circles $\sigma_{234}$, $\sigma_{134}$, $\sigma_{124}$, $\sigma_{123}$ all pass through one point $S_{1234}$.
\vskip 3pt

{\rm (2)} Let $\sigma_5$ be a fifth circle through $S$. Then the five points $S_{2345}$, $S_{1345}$, $S_{1245}$, $S_{1235}$, 
$S_{1234}$ all lie on one circle $\sigma_{12345}$.
\vskip 3pt

{\rm (3)} The six circles $\sigma_{23456}$, $\sigma_{13456}$, $\sigma_{12456}$, $\sigma_{12356}$, $\sigma_{12346}$, 
$\sigma_{12345}$ all pass through one point $S_{123456}$. 
\vskip 3pt

And so on.
\end{thm}

\noindent
We know the Levi graph of these configurations (Coxeter~\cite{Cox50, Cox61}).

\begin{lemma} \label{Levi-Clifford}
The Levi graph of the Clifford configuration of type $(2^{d-1}_d)$ is isomorphic to the $d$-cube graph.
\end{lemma}

\noindent
It turns out that our $V\!$-construction can be applied so as to obtain Clifford's configurations. 

\begin{thm}
The $V\!$-construction on a $d$-cube graph gives rise to an isometric  $(2^d_d)$ point-circle configuration in the plane. This configuration is 
disconnected and composed of two copies of isometric $(2^{d-1}_d)$ point-circle configurations which are isomorphic to a member of the same type
of  Clifford's infinite series of configurations.
\end{thm}

\noindent
\emph{Proof.} The $d$-cube graph is the Cartesian product of $d$ edge graphs $K_2$. According to~\cite{HP10}, it is a unit-distance graph. 
By Theorem \ref{thm:ud}, the $V\!$-construction applied on it gives rise to an isometric point-circle configuration $C$. Since the $d$-cube 
graph is bipartite, its Kronecker cover is composed of two disjoint isomorphic copies of the $d$-cube graph (by Proposition 1 in~\cite{IP08}). 
By Theorem~\ref{MainThm}, this Kronecker cover is the Levi graph of $C$. Since a configuration is uniquely determined by its Levi graph
(by Lemma~\ref{Levi}), it follows from Lemma~\ref{Levi-Clifford} that $C$ is in fact composed of two disjoint copies of Clifford configurations 
of type $(2^{d-1}_d)$.
\hfill $\square$
\bigskip

\begin{figure}[h!] 
  \begin{center}
  \subfigure[]{\hskip -4pt
  \includegraphics[width=0.375\textwidth]{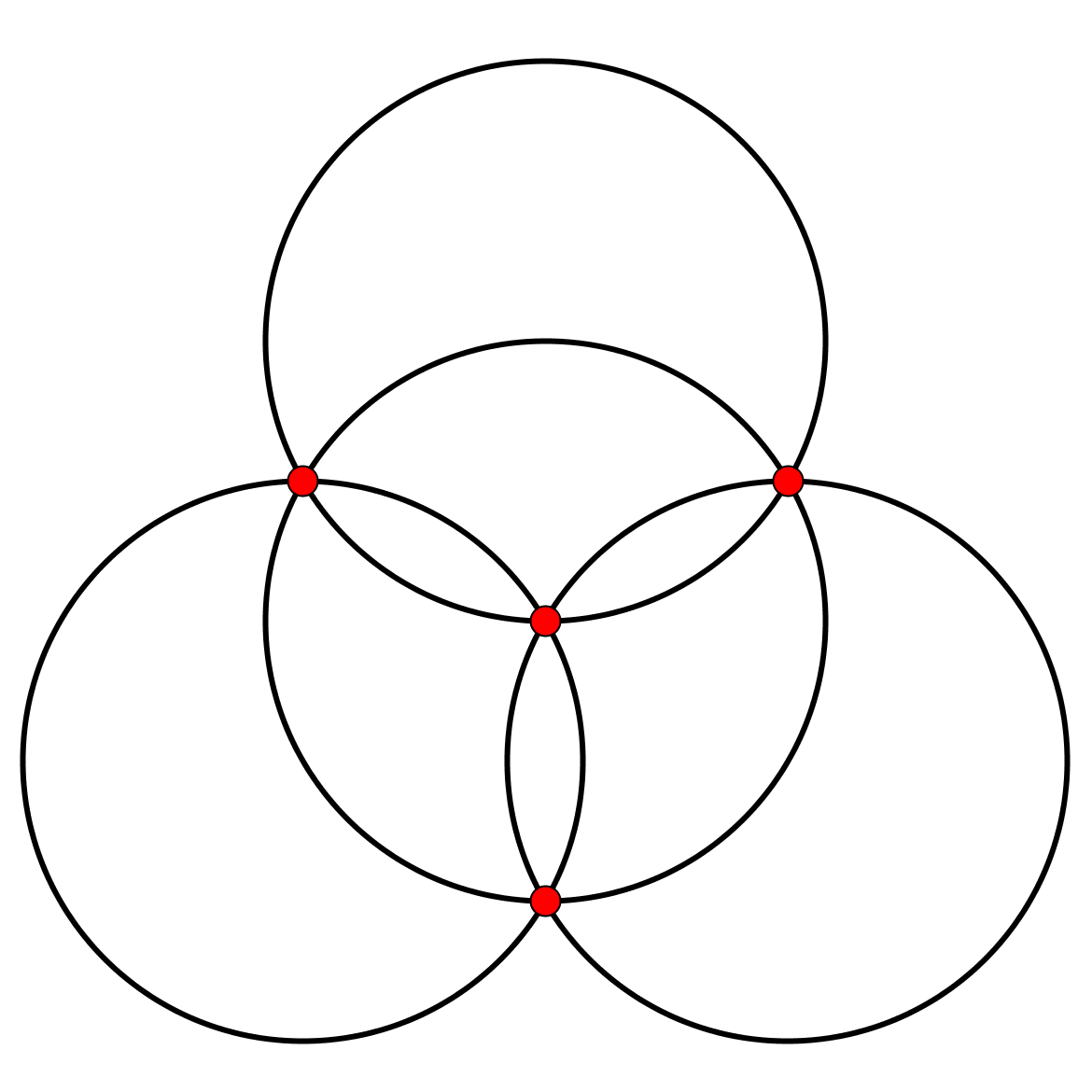}} \hskip 30pt
 \subfigure[]{\hskip -4pt
   \includegraphics[width=0.415\textwidth]{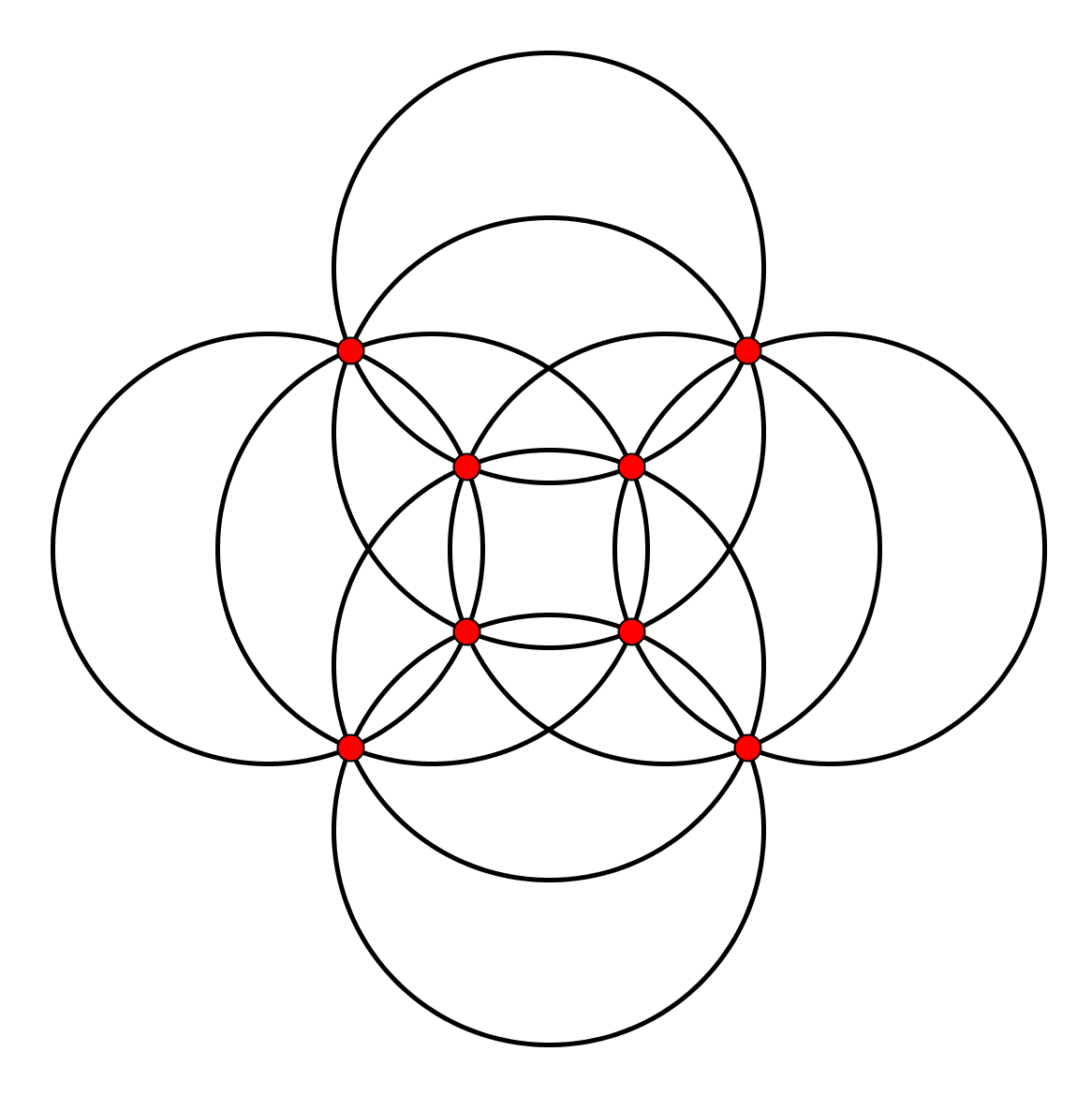}}
  \end{center}
  \caption{Isometric realization of Clifford configurations: $(4_3)$ derived from the 3-cube graph (a), and $(8_4)$ derived from the 4-cube graph (b).}
  \label{(4_3)-(8_4)}  
\end{figure}

\begin{figure}[h!] 
\begin{center}
\includegraphics[width=0.5\textwidth]{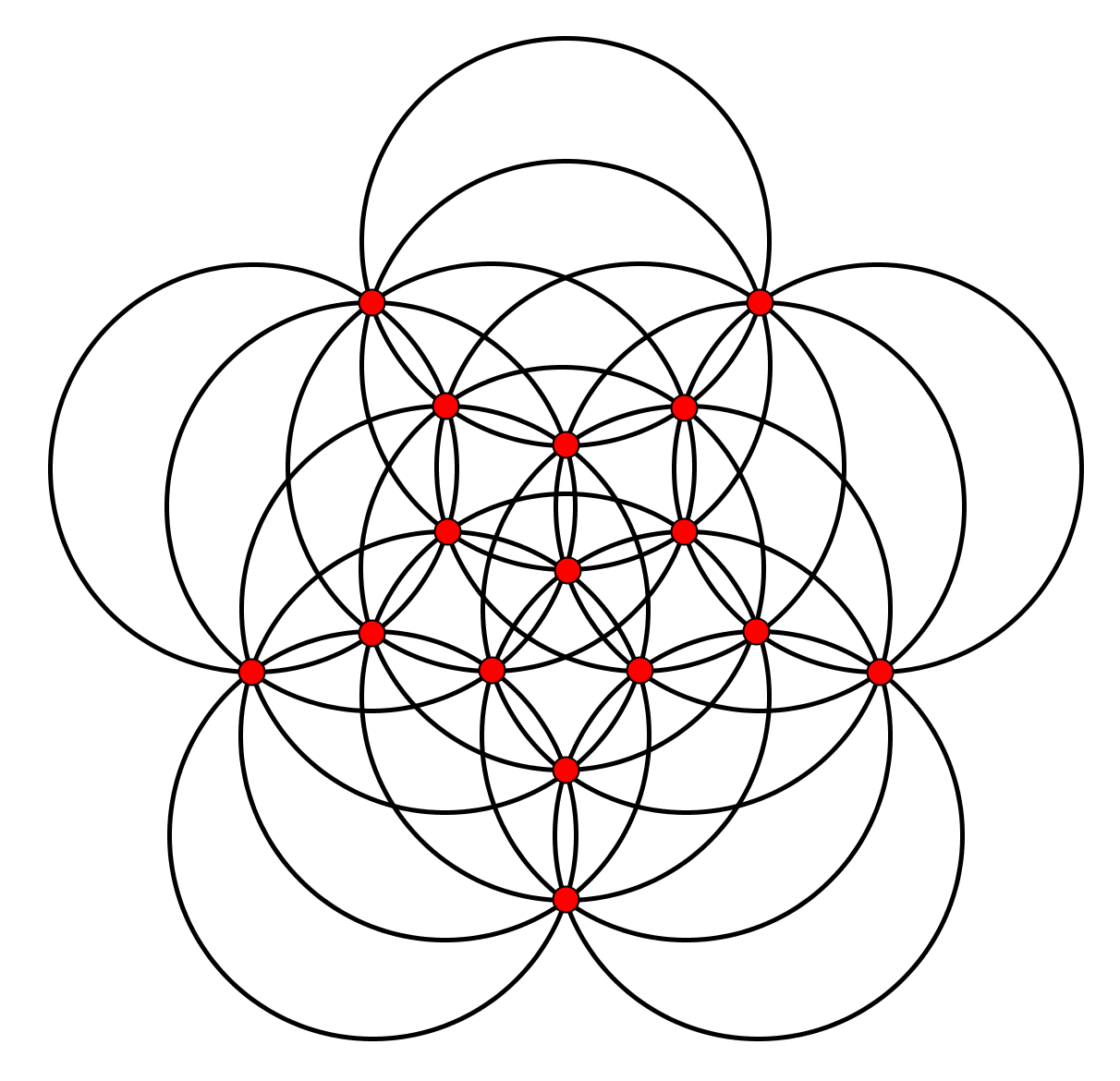}
\end{center}
 \caption{Isometric realization of Clifford's configuration $(16_5)$ derived from the 5-cube graph.}
  \label{(16_5)}  
\end{figure}

\noindent
Some smallest examples are depicted in Figures~\ref{(4_3)-(8_4)} and~\ref{(16_5)}.

It is easy to see that the $d$-cube graph can be realized in the plane as a unit-distance graph in continuum many ways. In fact, take an arbitrary vertex 
and place it in the center of a unit circle. Its first-neighbours can be placed in different positions on the circle. Positions of the remaining vertices are 
then uniquely determined by sequences of rhombuses. This immediately gives the following corollary.    

\begin{cor}\label{IsometricRealization}
Every Clifford configuration is realizable as a movable isometric point-circle configuration.
\end{cor}

\begin{remark}
Realizability of Clifford's configurations with circles of equal size is already known from~\cite{Zie40} (see also~\cite{Bab}). 
Our approach provides an independent proof of  this result.
\end{remark}

\section{Three new infinite classes of point-circle configurations} \label{KneserConfig}

\begin{figure}[h!] 
  \begin{center}
  \subfigure[]{\hskip -4pt
  \includegraphics[width=0.375\textwidth]{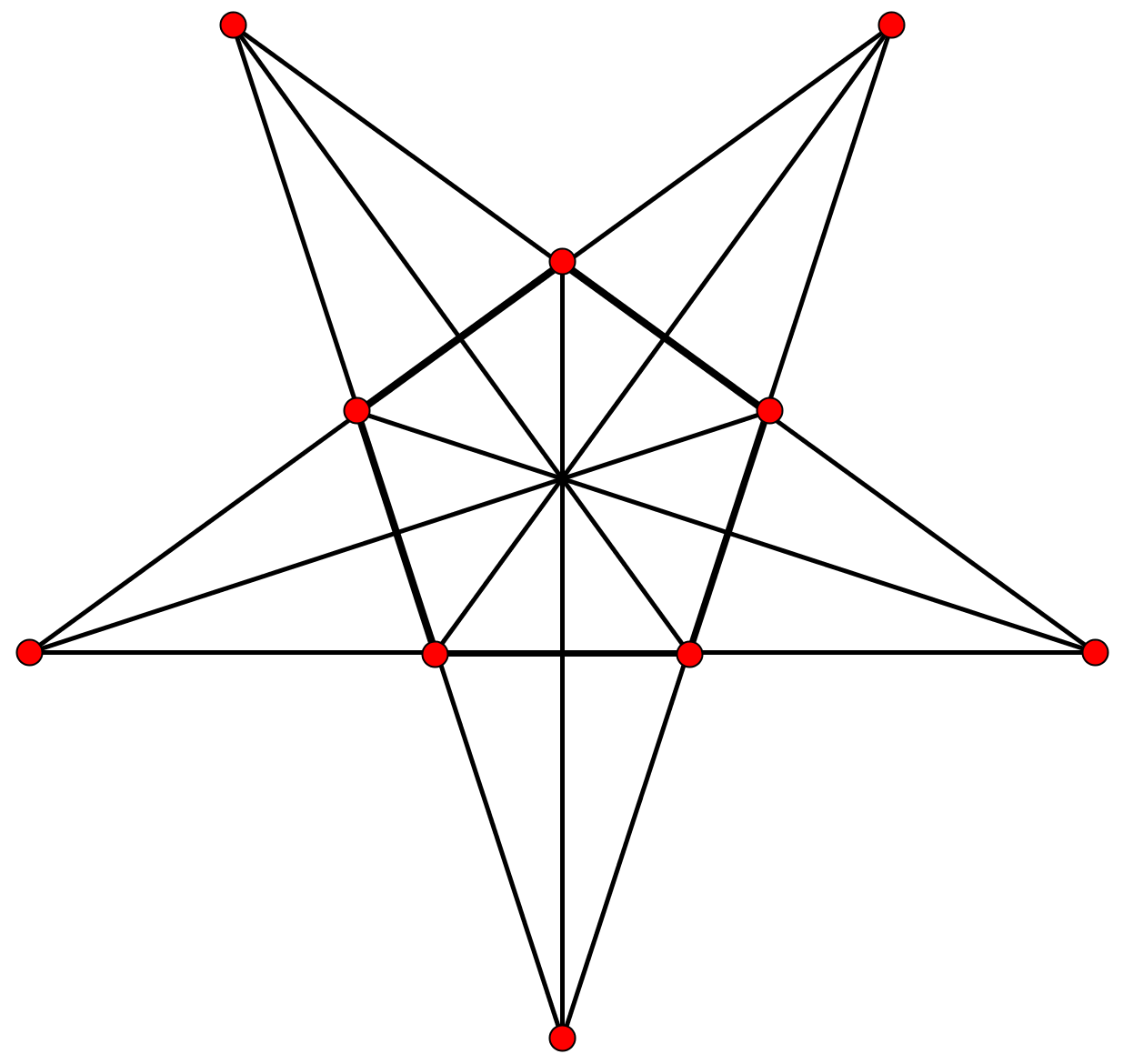}} \hskip 30pt
 \subfigure[]{\hskip -4pt
   \includegraphics[width=0.45\textwidth]{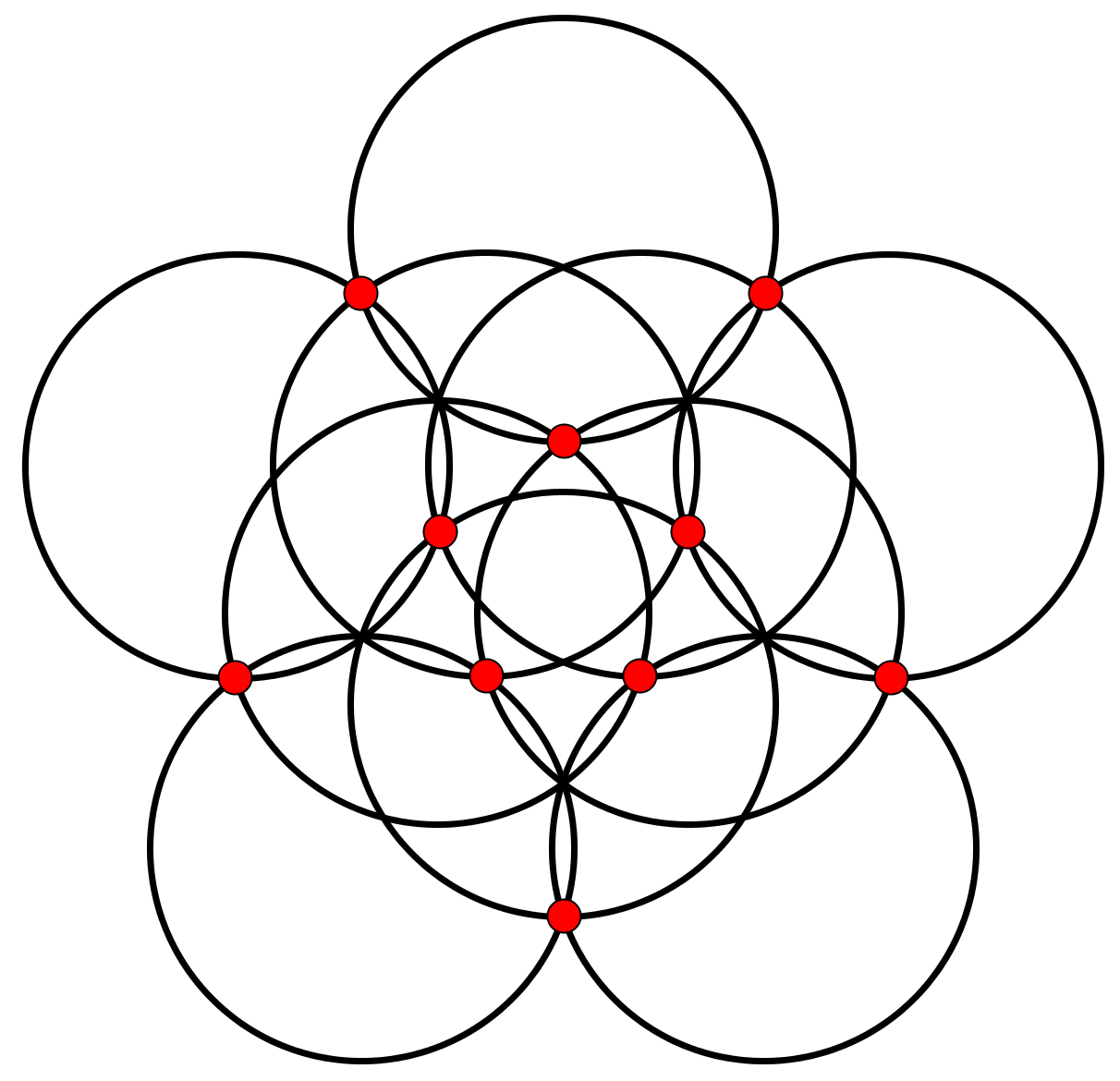}}
  \end{center}
  \caption{Isometric point-circle realization of Desargues' configuration arising from a non-unit-distance representation of the Petersen graph.
                (Note that this is a degenerate representation in the sense that there are edges overlapping along the sides of the inner pentagon.)}
  \label{NUD-Desargues}  
\end{figure}

\noindent
We start from the following observation. When applying the $V\!$-construction so as to obtain an isometric point-circle realization 
of Desargues' configuration, the underlying Petersen graph need not be represented in unit-distance form. Indeed, compare our 
Figures~\ref{UD-Desargues} and~\ref{NUD-Desargues}. (We remark that the version depicted in Figure~\ref{NUD-Desargues}b is 
precisely the same as one of the components of the $(20_3)$ configuration in Figure~\ref{DoubleDesargues}b.)

\begin{figure}[h!] 
  \begin{center} 
  \includegraphics[width=0.55\textwidth]{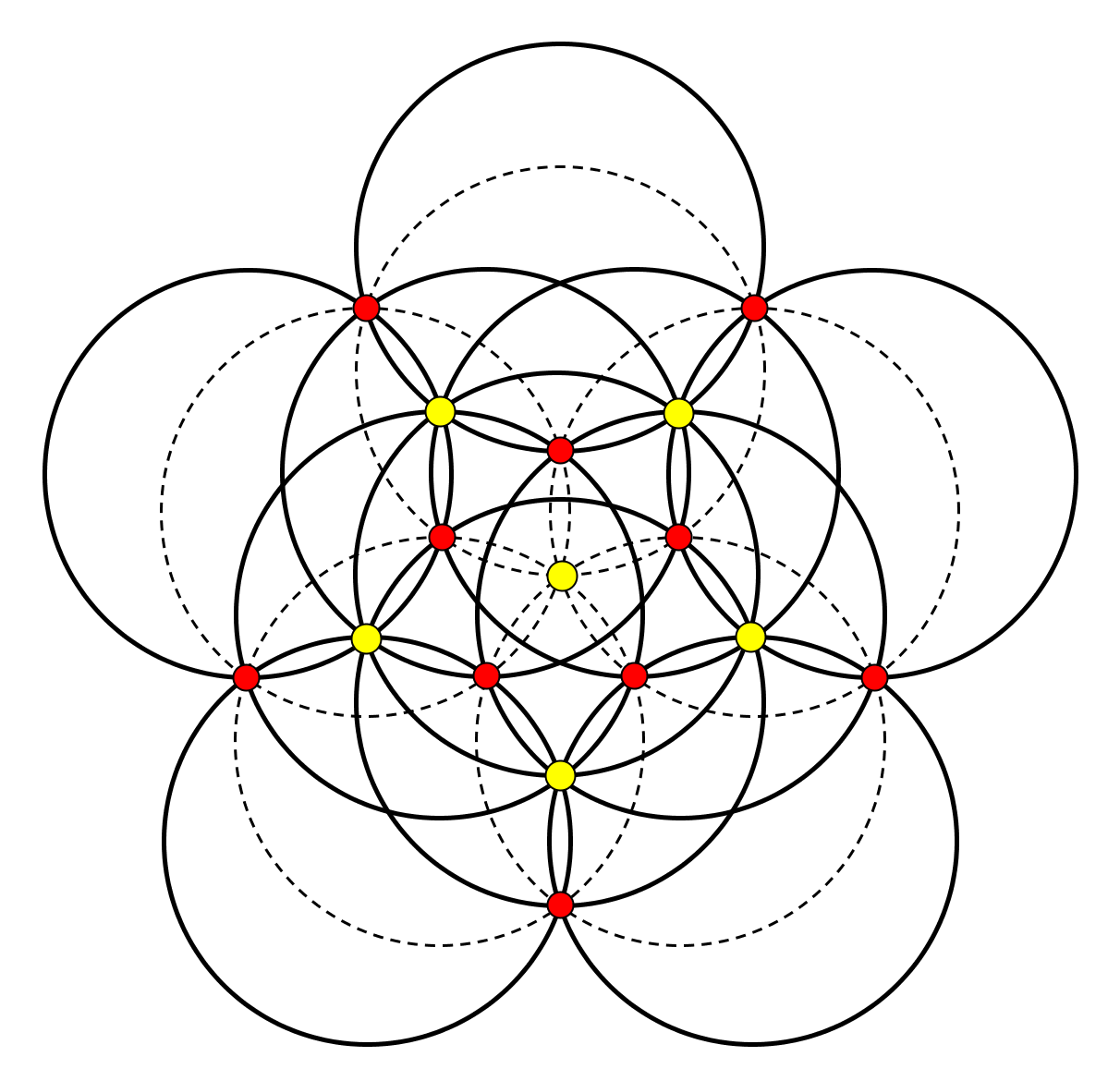}
  \end{center}
  \caption{The Desargues point-circle configuration extended to the Clifford configuration $(16_5)$.}
  \label{DesarguesExtended}  
\end{figure}

Now Figure~\ref{NUD-Desargues}b suggests that this latter realization can be extended to a Clifford configuration of type $(16_5)$. 
Figure~\ref{DesarguesExtended} shows that such an extension is in fact possible (see also Figure~\ref {(16_5)}).
It turns out that this is a particular case of a more general relationship. 

Before formulating it, recall that the \emph{Kneser graph} $K(n, k)$ has as vertices the $k$-subsets of an $n$-element set, where two vertices 
are adjacent if the $k$-subsets are disjoint~\cite{GR01}. The Kneser graph $K(2n -1,n - 1)$ is called an \emph{odd graph} and is denoted by 
$O_n$. In particular, $O_3 = K(5,2)$ is isomorphic to the Petersen graph. The \emph{bipartite Kneser graph} $H(n, k)$ has as its bipartition 
sets the $k$- and $(n-k)$-subsets of an $n$-element set, respectively, and the adjacency is given by containment. Although the following relationship 
is well-known, we give a short proof of it.

\begin{lemma} \label{KneserCover}
The bipartite Kneser graph $H(n, k)$ is the Kronecker cover of the Kneser graph $K(n, k)$.
\end{lemma}

\noindent
\emph{Proof.}
Let $A$ and $B$ be two $k$-subsets and let $A'$ and $B'$ be their respective $(n-k)$-complements. Clearly
$A$ is adjacent to $B$ in $K(n,k)$ if and only if $A$ is adjacent to $B'$ and $B$ is adjacent to $A'$ in $H(n,k)$,
 and the result follows readily.
\hfill $\square$
\vskip 8pt

The bipartite Kneser graph $H(2n-1, n-1)$ is also known as the \emph{revolving door graph}, or \emph{middle-levels graph}; the latter name 
comes from the fact that it is a special subgraph of the $(2n-1)$-cube graph $Q_{2n-1}$ (considering $Q_{2n-1}$ as the Hasse diagram of the 
corresponding Boolean lattice)~\cite{QB, SS}. It is a regular graph with degree $n$. Note that middle-levels graph is called a \emph{medial 
layer graph} in~\cite{MPSI} and is defined for any abstract polytope of odd rank.

\begin{thm}
For all $n\ge3$, there exists an isometric point-circle configuration of type
$$
\left(
\dbinom{2n-1}{n-1}
_{n}\,
\right).
$$
It is a subconfiguration of the Clifford configuration of type $\left(2^{2n-2}_{2n-1}\right)$. It can be obtained from the odd graph $O_n$ by
$V\!$-construction.
\end{thm}

\noindent
\emph{Proof.}  Let $C$ be an incidence structure obtained from the odd graph $O_n$ by $V\!$-construction. By Theorem~\ref{MainThm}, the 
Levi graph of $C$ is the Kronecker cover of  $O_n$. Lemma~\ref{KneserCover} implies that it is the bipartite Kneser graph $H(2n-1, n-1)$. 
Since this graph is a subgraph of the $(2n-1)$-cube graph $Q_{2n-1}$, from Lemma~\ref{Levi} follows that $C$ is isomorphic to a subconfiguration 
of the Clifford configuration of type $\left(2^{2n-2}_{2n-1}\right)$. Hence it can be realized as a planar point-circle configuration. The type of 
this configuration follows from the definition of $O_n$. Furthermore, Corollary~\ref{IsometricRealization} implies that this configuration also 
has an isometric realization. 
\hfill $\square$
\vskip 8pt

\begin{figure}[h!] 
  \begin{center} 
  \includegraphics[width=0.75\textwidth]{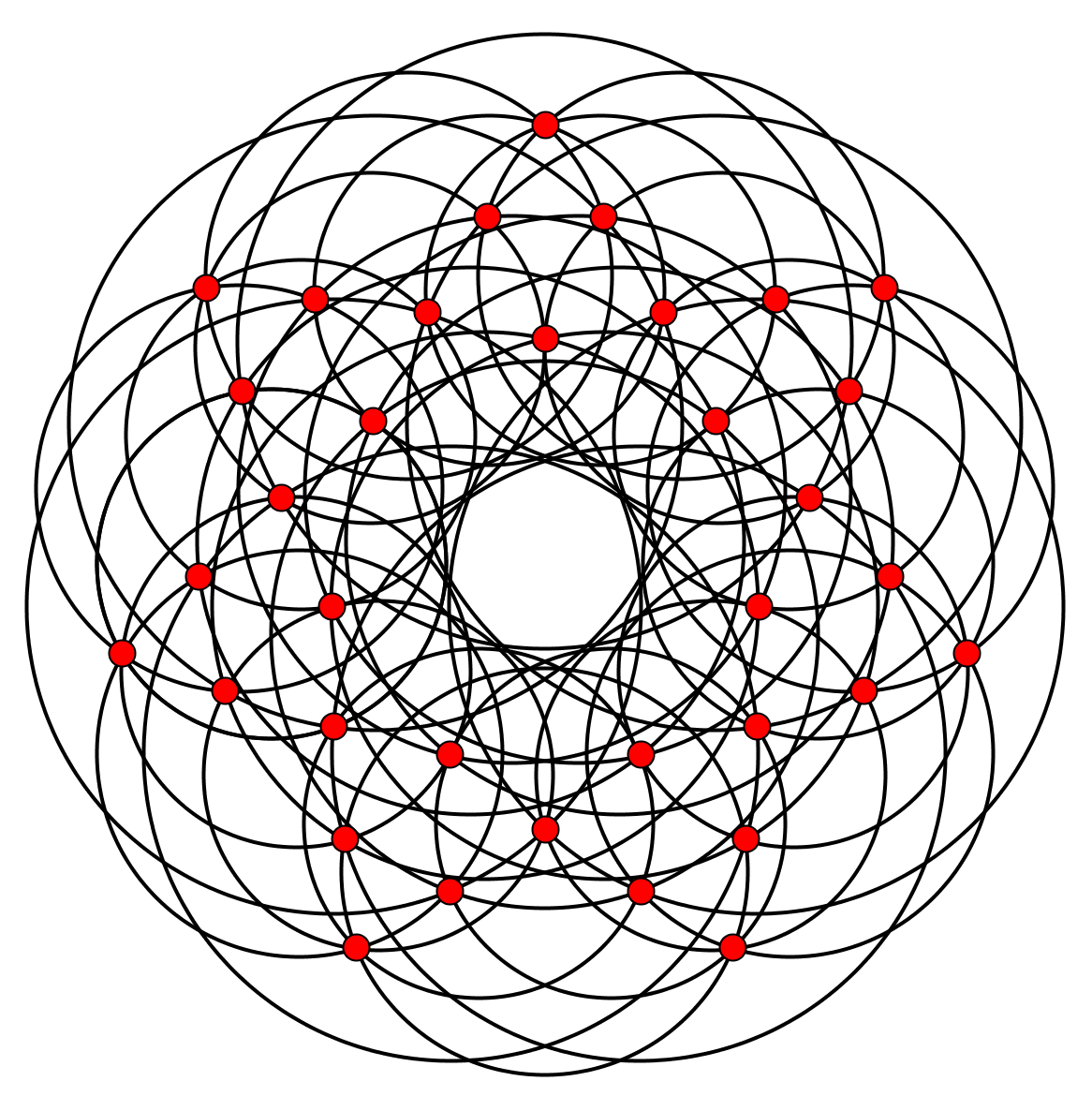}
  \end{center}
  \caption{A point-circle realization of Danzer's $(35_4)$ configuration.} \label{Danzer}  
\end{figure}

In the particular case of $n=4$ we have $(35_4)$, which provides a point-circle realization of Danzer's $(35_4)$ point-line configuration 
(see Figure~\ref{Danzer} for a non-isometric version). On this latter, Gr\"unbaum wrote in 2008~\cite{Gru08}: ``It seems that any representation 
of Danzer's configuration is of necessity so cluttered and unhelpful for visualization that no attempt to present it has ever been made." (see 
also~\cite{GR90}). 
We emphasize the geometric symmetry of this realization, which is the highest possible in the planar case; namely, $D_7$. 

Our next new class also consists of isometric point-circle configurations.

\begin{thm}
For any $N$ and any $d > 2$ there exists an isometric  $(n_d)$ point-circle configuration with $n > N$.
\end{thm}

\noindent
\emph{Proof.}  Take the Cartesian product of a long odd cycle $C_N$ and a $(d-2)$-dimensional cube graph. This is a unit-distance graph.  
Apply the $V\!$-con\-struc\-tion to it. \hfill $\square$
\vskip 8pt

Finally, we construct an infinite series of non-isometric  $(n_4)$ point-circle configurations. We start from a prism $P$ over an $n$-gon 
($n \ge 3$) (the corresponding graph is also called a \emph{circular ladder}). Then we take its \emph{medial} $Me(P)$~\cite{GR01,
PR,FP}, i.e.\ a new polyhedron such that its vertices are the midpoints of the original edges, and for each original vertex, the midpoints 
of the edges emanating from it are connected by new edges, forming a 3-cycle. In terms of solid geometry, the medial corresponds to 
a truncation of a right $n$-sided prism such that each truncating plane at a vertex is spanned by the midpoints of the edges incident with 
the given vertex (``deep vertex truncation", see~\cite{Zie07}). Note that in the particular case when the prism is the cube, its medial 
is the Archimedean solid called a \emph{cuboctahedron}. Accordingly, we define the \emph{generalized cuboctahedron graph} as the 
1-skeleton of $Me(P)$, and denote it by $CO(n)$. Note that this graph can equivalently be defined as the line graph of the prism graph.

Observe that $CO(n)$ is an $4$-valent regular graph with $3n$ vertices. Moreover, it has a representation in the plane such that it exhibits 
the symmetry of a regular $n$-gon (thus its symmetry group is $D_n$); in this case its vertices lie on three concentric circles, $n$ vertices 
on each. It follows that the first-neighbour sets of the vertices are concyclic, hence the $V\!$-construction can be applied. Thus we obtain the 
following result.

\begin{thm} \label{GenCub}
For any $n\ge3$, there exists a $((3n)_4)$ point-circle configuration obtained from the generalized cuboctahedron graph 
$CO(n)$ by $V\!$-construction. It can be realized in the plane so that its symmetry group is the dihedral group $D_n$.
\end{thm}

\begin{figure}[h!] 
  \begin{center}
  \subfigure[]{\hskip -4pt
  \includegraphics[width=0.475\textwidth]{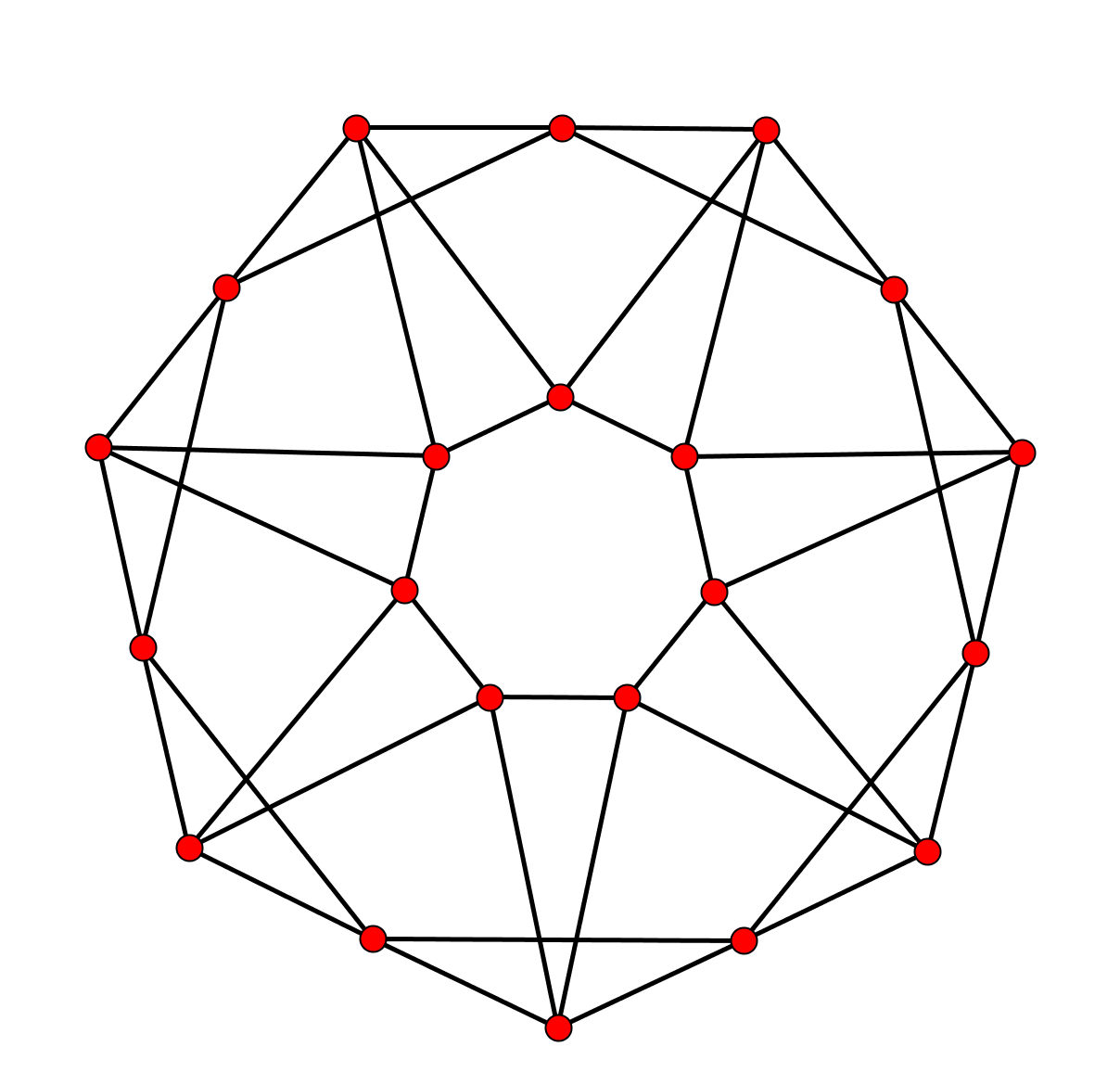}} \hskip 15pt
 \subfigure[]{\hskip -4pt
   \includegraphics[width=0.465\textwidth]{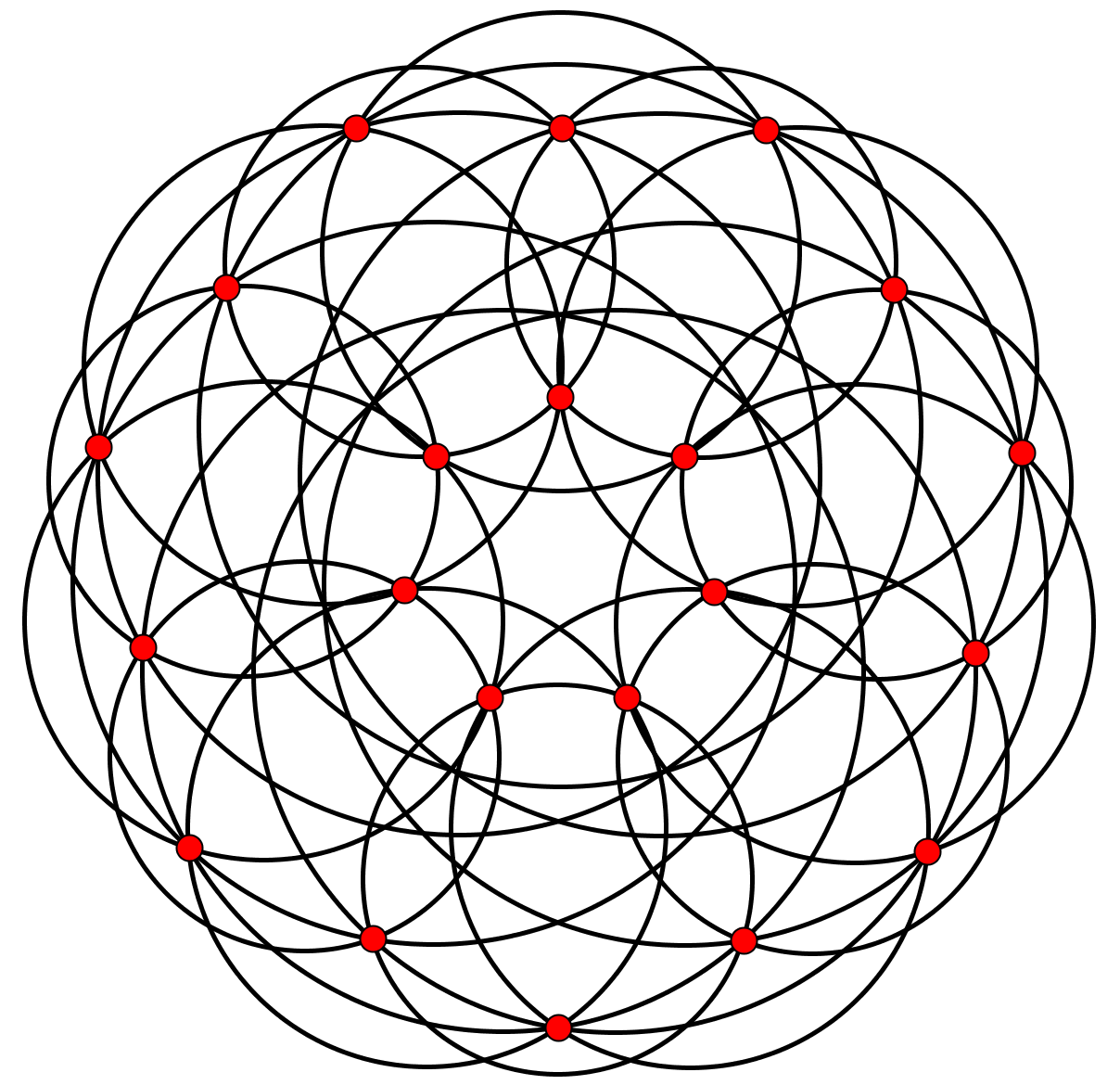}}
  \end{center}
  \caption{A representation of the generalized cuboctahedron graph with $D_7$ symmetry (a), 
                and a $(21_4)$ point-circle configuration obtained from it by $V\!$-construction (b).}
  \label{Cuboctahedral}  
\end{figure}

\noindent
An example with $n=7$ is depicted in Figure~\ref{Cuboctahedral}. It is an open question if any member of the infinite 
series of these $((3n)_4)$ configurations has an isometric realization.

\begin{figure}[h!] 
  \begin{center} \hskip 8pt
  \subfigure[]{\hskip- 4pt
  \includegraphics[width=0.425\textwidth]{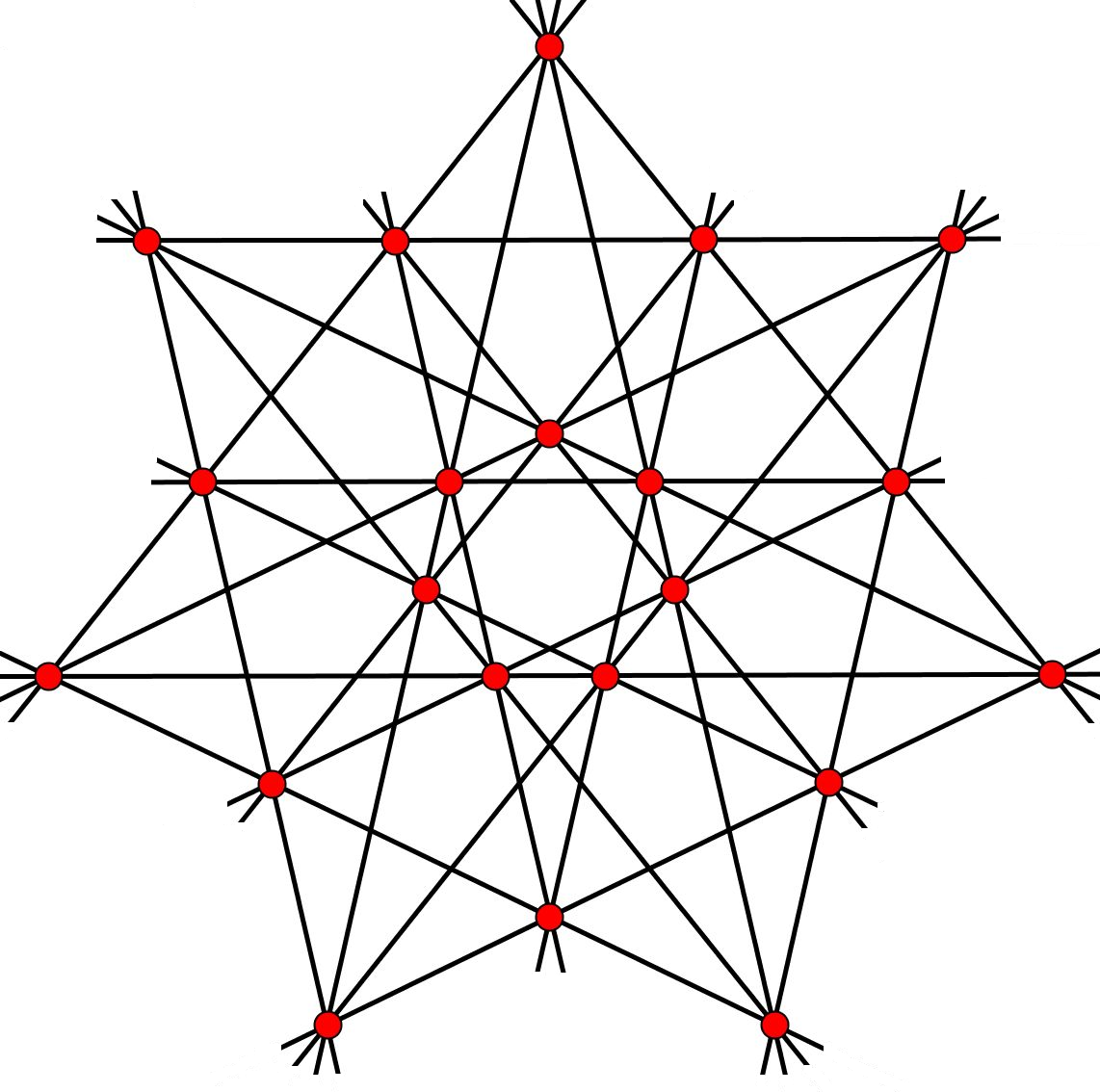}} \hskip 18pt
 \subfigure[]{\hskip -4pt
   \includegraphics[width=0.5\textwidth]{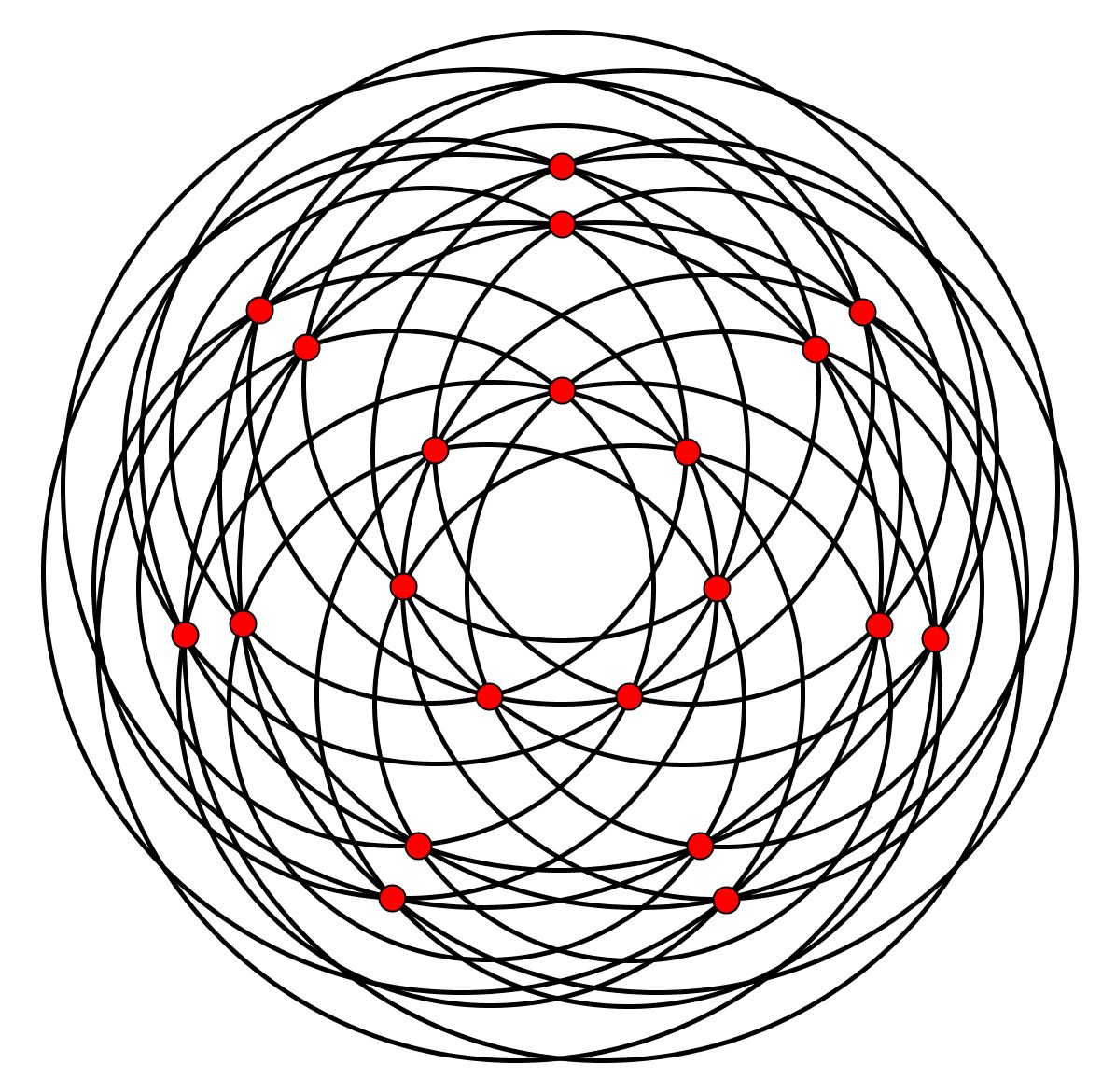}}
  \end{center}
  \caption{The point-line $(21_4)$ configuration of Gr\"unbaum and Rigby (a), 
                one of its point-circle realizations with dihedral symmetry (b).}
  \label{RigbyGrunbaum}  
\end{figure}

On the other hand, the mutual position of the points on the three orbits makes it possible to arrange the circles in several different ways, so as to 
obtain new, pairwise non-isomorphic $((3n)_4)$ configurations. Here we do not investigate this possibility in detail. Instead, we just present an 
example, also of type $(21_4)$ (non-isomorphic with the previous one), whose original point-line version is remarkable for several reasons (see 
Figure~\ref{RigbyGrunbaum}). We only mention here that it goes back to Felix Klein, 1879 (for further details, see~\cite{GR90}); on the other 
hand, its first graphic depiction only appeared in 1990~\cite{Gru08, GR90}. This configuration also motivated the authors of \cite{MP} to 
present some geometric representations of a certain family of configurations that became later known as polycyclic configurations~\cite{BP}.

We note that several other already known families of graphs can serve as basis for obtaining new point-circle configurations by V-construction;
just to mention some of them: generalized Petersen graphs~\cite{PR}, $I$-graphs~\cite{BPZ}. In addition, the 1-skeleton of equivelar polyhedra 
is also a regular graph (see e.g.~\cite{GW}); hence, finding suitable planar representations among them may also be promising in this respect.

\section{Some comparisons beween different realizations of con\-figur\-a\-tions} \label{Comparisons}

\noindent
Comparing point-line and point-circle configurations, several questions arise, in particular, when different kinds of geometric realization of 
the same abstract configuration is considered. First we make the following conceptual distinction. Clearly, every point-line configuration can 
be transformed into a point-circle configuration by some suitable inversion. However, in this case, all the circles will have a common point
(the inverse image of the point at infinity). To rule out this case, we use the term {\em improper point-circle configuration\/}. Accordingly, 
we call a point-circle configuration {\em proper\/} if its circles are not all incident with a common point. Clearly, all our examples presented 
in the previous sections are proper point-circle configurations. In what follows, we shall also speak about such configurations, and mostly omit 
the attribute ``proper".

A simple consequence of Proposition~\ref{n_3} is that by suitable displacement of the points, any planar ${n_3}$ point-line configuration can 
transformed into a point-circle configuration. For an incidence number larger than 3, it is more difficult to decide the existence of a point-circle 
representation of a point-line configuration.

\begin{problem}
For $k \ge 4$, find an $(n_k)$ point-line configuration which cannot be represented by a proper point-circle configuration.
\end{problem}

The converse problem, in general, can also be quite difficult. However, here we know several examples. One of the oldest one is Miquel's
$(8_3, 6_4)$ (for a simple proof why it has no point-line representation, see~\cite{PS}). The infinite series of Clifford configurations also 
provides quite old, and balanced examples. In fact, since all the higher members contain, as a subconfiguration, the initial member of type
$(4_3)$, they cannot be represented by point-line configurations.

In the particular case of incidence number $k=4$, we have the following lower bound (a result of Bokowski and Schewe~\cite{BS}).

\begin{thm}
For $n\le 17$, there are no geometric point-line configurations $(n_4)$.
\end{thm}

As a consequence, consider e.g.\ the generalized Petersen graph $GP(n,r)$~\cite{PR}. For $n\le 8$ it yields, by $V\!$-construction, a 
point-circle configuration which has no point-line representation.

In Section~\ref{PointCircle} we introduced the notion of an isometric point-circle configuration. We may impose two further conditions, which,
together with the former, determine a particularly nice class of configurations. We call a point-circle configuration $\mathcal C$ {\em lineal\/} if 
two circles meet in at most one configuration points. Furthermore,  $\mathcal C$ is called {\em determining\/} if the set of points of  $\mathcal C$
coincides with the set of points in which more than two circles of $\mathcal C$ meet.

Note that these two conditions differ in the sense that the former determines a property on more abstract level, i.e.\ $\mathcal C$ is lineal if and 
only if it is isomorphic to a combinatorial configuration which is likewise lineal (we may call such a property of a geometric configuration {\em 
intrinsic}). On the other hand, the latter may depend on a particular representation of $\mathcal C$. For example, Figure~\ref{UD-Desargues}b 
shows a determining representation of the Desargues configuration, while that in Figure~\ref{NUD-Desargues}b is non-determining. (Such a 
property may be called {\em extrinsic\/}; note that being isometric is another example of an extrinsic property in this sense.)

Now we call $\mathcal C$ {\em perfect\/} if it is lineal, isometric and determining. For example, the Desargues configuration in 
Figure~\ref{UD-Desargues}b is perfect.  A similar question can be posed for point-line configurations.

\begin{problem}
Which configurations of points and lines can be realized as perfect point-circle configurations?
\end{problem}

\begin{figure}[h!] 
  \begin{center}
  \subfigure[]{\hskip -4pt
  \includegraphics[width=0.425\textwidth]{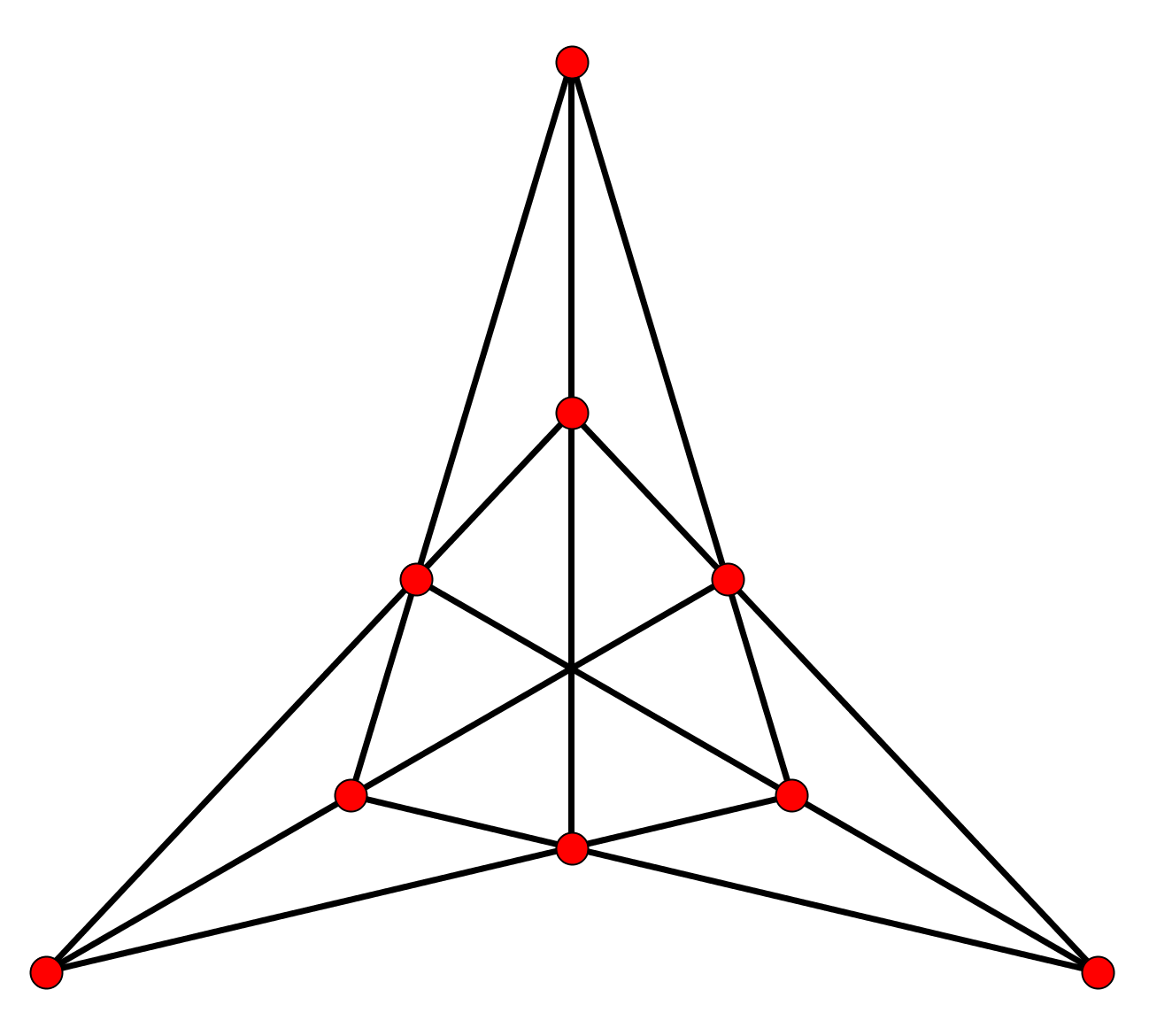}} \hskip 20pt
 \subfigure[]{\hskip -4pt
   \includegraphics[width=0.4\textwidth]{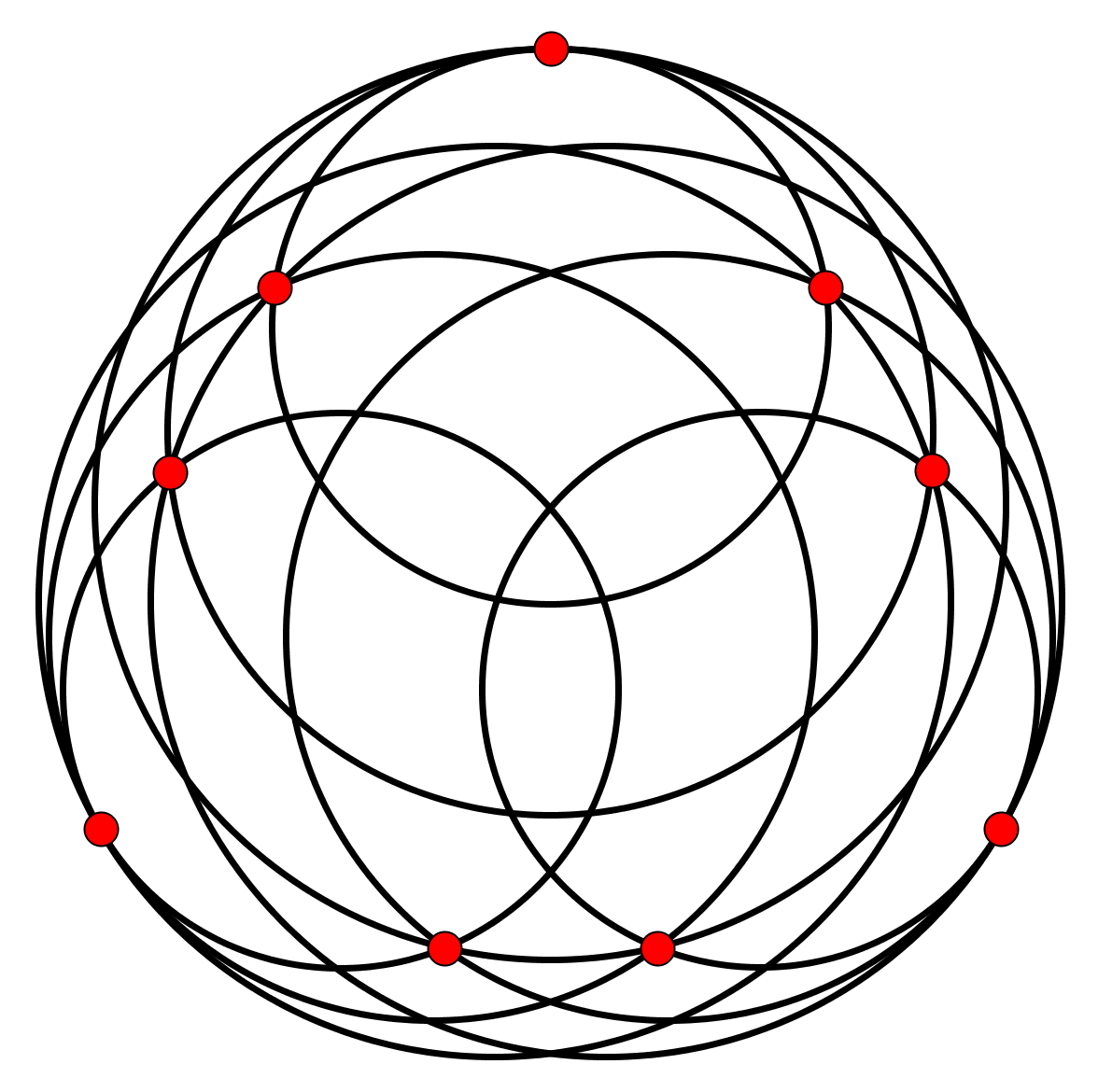}}
  \end{center}
  \caption{Two different realizations of Pappus' configuration with symmetry $D_3$.}
  \label{Pappus}  
\end{figure}

Geometric symmetry is also an interesting property which is worth investigating when different realizations of the same abstract configuration
are compared. Are all symmetries of a point-line configuration realizable in its representation by points and circles? Of course, the converse 
question can also arise. Here we only mention that e.g.\ for the Pappus configuration not only its realization by lines can exhibit the maximal 
possible symmetry ($D_3$), but it can also be realized by circles with the same symmetry (see Figure~\ref{Pappus}).

On the other hand, it is a remarkable fact that while the Desargues configuration can be represented by points and lines with symmetry group 
either $C_5$ or $D_5$ (see Figures~\ref{UD-Desargues}b and~\ref{NUD-Desargues}b, respectively), its classical point-line version can 
exist with neither of these symmetries. This follows from the theory developed in the paper~\cite{BPZ} on $I$-graphs and the corresponding 
configurations. (We note that geometric realization of certain combinatorial objects with maximal symmetry is, in general, a problem which is 
far from trivial, see e.g.~\cite{GG09a}, and the references therein.)

\begin{figure}[h!] 
  \begin{center}
  \subfigure[]{\hskip -4pt
  \includegraphics[width=0.435\textwidth]{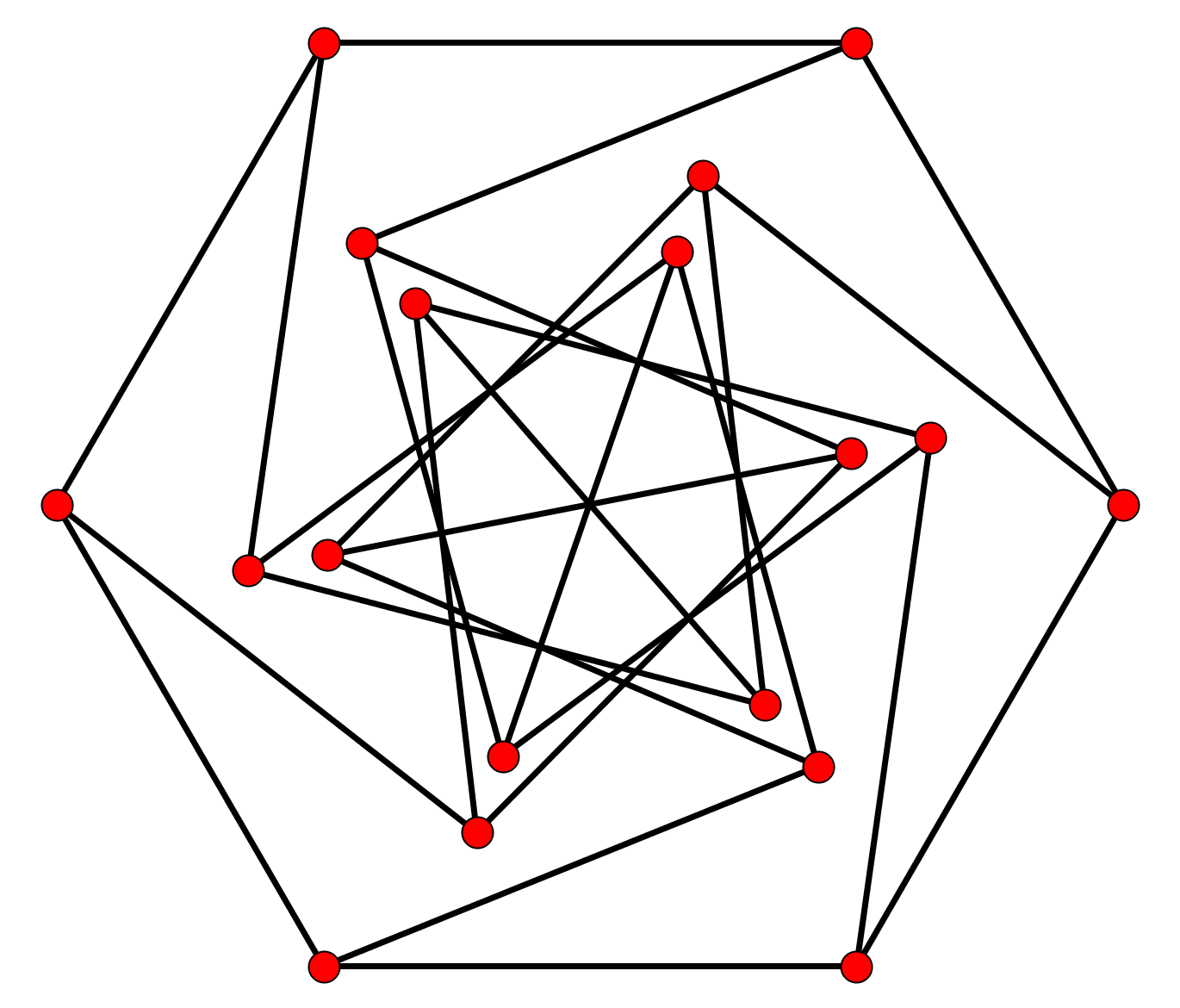}} \hskip 20pt
 \subfigure[]{\hskip -4pt
   \includegraphics[width=0.4\textwidth]{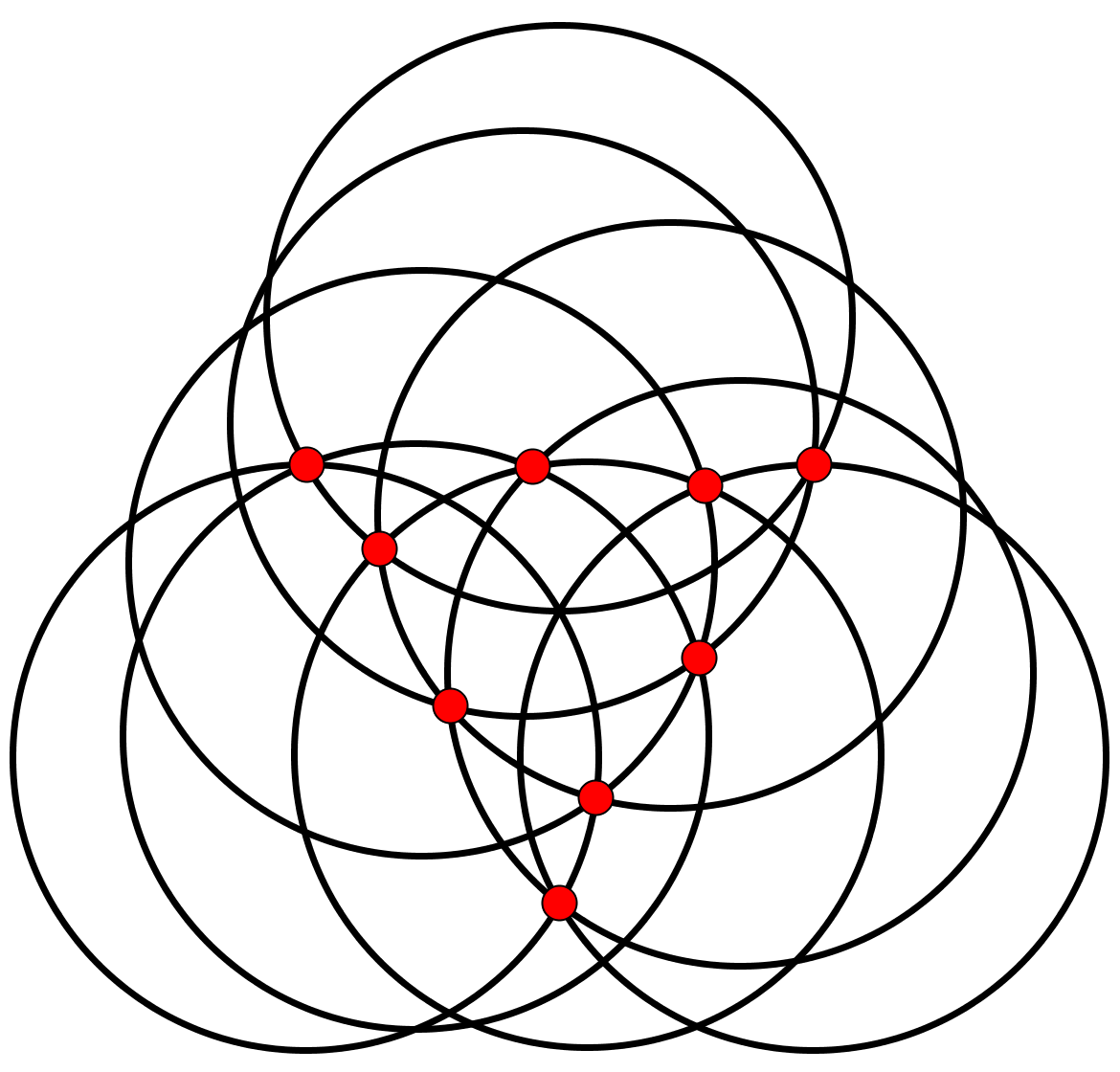}}
  \end{center}
  \caption{Unit-distance representation of the Pappus graph (a), and isometric point-circle representation of the Pappus configuration derived from it (b).}
  \label{PappusIsometric}  
\end{figure}
Considering point-circle realizations of the two oldest configurations, yet another difference occurs. Note that while Desargues' configuration 
has a perfect realization (shown by Figure~\ref{UD-Desargues}b), the realization of Pappus' configuration shown by Figure~\ref{Pappus}b 
is not isometric (thus it is not perfect). On the other hand, when constructing an isometric representation, we find that it loses the property being 
determining (see the central triple crossing point in Figure~\ref{PappusIsometric}b). This version is obtained from a unit-distance representation 
of the Pappus graph (see Figure~\ref{PappusIsometric}a), using Corollary~\ref{TwoCopies}. Note that the symmetry reduces here to $C_3$
(for a representation of the Pappus graph with $D_3$ symmetry see e.g.~\cite{PR}, Figure 21).

\section*{Acknowledgement}
The authors would like to thank Marko Boben for careful reading of the man\-u\-script and for pointing out that some circles in the $(21_4)$ 
configuration on Figure~\ref{Cuboctahedral}(b) meet in two points, which is an easy proof that the configuration is not isomorphic to Klein's 
$(21_4)$ configuration from \cite{GR90}; and last, but not least, for checking that Danzer's combinatorial $(35_4)$ configuration actually 
admits a geometric realization and is hidden in the renowned Pascal's \emph {Hexagrammum mysticum}~\cite{Klu}.

\end{document}